\documentclass[a4paper,12pt]{article}

\usepackage{amscd}
\usepackage{amsfonts}
\usepackage{amsmath}
\usepackage{amssymb}
\usepackage{amsthm}
\usepackage{ascmac}
\usepackage{here}
\usepackage{makeidx}
\usepackage{mathrsfs}
\usepackage{slashbox}
\usepackage{txfonts}

\allowdisplaybreaks

\newcommand{\C}{\mathbb{C}}
\newcommand{\F}{\mathbb{F}}
\newcommand{\N}{\mathbb{N}}
\newcommand{\Q}{\mathbb{Q}}
\newcommand{\R}{\mathbb{R}}
\newcommand{\Z}{\mathbb{Z}}

\newcommand{\A}{\mathscr{A}}
\newcommand{\B}{\mathscr{B}}

\newcommand{\M}{\mathscr{M}}

\theoremstyle{plain}
\newtheorem{thm}{Theorem}[section]
\newtheorem{lmm}[thm]{Lemma}
\newtheorem{prp}[thm]{Proposition}
\newtheorem{crl}[thm]{Corollary}
\theoremstyle{definition}
\newtheorem{dfn}[thm]{Definition}

\newtheorem{exm}[thm]{Example}

\def\Set#1#2{\left\{ #1 \ \middle| \ #2 \right\}}
\def\set#1#2{\{ #1 \mid #2 \}}
\def\m#1{\text{\rm #1}}

\title{
Characterisation of the Berkovich Spectrum of the Banach Algebra of Bounded Continuous Functions
}
\author{Tomoki Mihara}
\date{}

\begin{document}

\maketitle
\begin{abstract}
For a complete valuation field $k$ and a topological space $X$, we prove the universality of the underlying topological space of the Berkovich spectrum of the Banach $k$-algebra $\m{C}_{\m{bd}}(X,k)$ of bounded continuous $k$-valued functions on $X$. This result yields three applications: a partial solution to an analogue of Kaplansky conjecture for the automatic continuity problem over a local field, comparison of two ground field extensions of $\m{C}_{\m{bd}}(X,k)$, and non-Archimedean Gel'fand theory.
\end{abstract}

\tableofcontents
\newpage

\section{Introduction}
\label{Introduction}

A non-Archimedean analytic space plays an important role in various studies in modern number theory. There are several ways to formulate a non-Archimedean analytic space, and one of them is given by Berkovich in \cite{Ber1} and \cite{Ber2}. Berkovich introduced the spectrum $\M_k(\A)$ of a Banach algebra $\A$ over a complete valuation field $k$. The space $\M_k(\A)$ is to a Banach algebra $\A$ what $\m{Spec}(A)$ is to a ring $A$. We note that $\M_k(\A)$ is called the Berkovich spectrum in modern number theory, but the same notion is originally defined by Bernard Guennebaud in \cite{Gue}. The class of Banach algebras topologically of finite type over a complete valuation field is significant in analytic geometry, just as the class of algebras of finite type over a field is significant in algebraic geometry. A Banach algebra topologically of finite type is called an affinoid algebra, and the Berkovich spectrum of an affinoid algebra is called an affinoid space. The space $\M_k(\A)$ is a compact Hausdorff G-topological space. For the notion of G-topology, see \cite{BGR}. Berkovich formulated an analytic space by gluing affinoid spaces with respect to a certain G-topology, just as Grothendieck did a scheme by gluing affine schemes with respect to the Zariski topology. We remark that an affinoid space is studied well, while few properties are known for the Berkovich spectrum of a general Banach algebra.

\vspace{0.1in}
Throughout this paper, $X$ and $k$ denote a topological space and a complete valuation field respectively. Here a {\it valuation field} means a field endowed with a valuation of height at most $1$, and we allow the case where the valuation is trivial. We study the underlying topological space of the Berkovich spectrum $\m{BSC}_k(X)$ of the Banach algebra $\m{C}_{\m{bd}}(X,k)$ of bounded continuous $k$-valued functions on $X$. In Theorem \ref{main theorem}, we prove that $\m{BSC}_k(X)$ is naturally homeomorphic to the Stone space $\m{UF}(X)$ associated to $X$, where $\m{UF}(X)$ is a topological space under $X$ (Definition \ref{under}) constructed using the set of ultrafilters of a Boolean algebra associated to $X$. This homeomorphism is significant because $\m{UF}(X)$ is an initial object in the category of totally disconnected compact Hausdorff spaces under $X$ (Definition \ref{initial}). As a consequence, $\m{BSC}_k(X)$ satisfies the same universality, and hence is independent of $k$. We note that Banaschewski proved the existence of such an initial object only for zero-dimensional spaces in \cite{Ban} Satz 2, while we deal with a general topological space in this paper. We also remark that many of our results are verified by Alain Escassut and Nicolas Ma\"inetti in \cite{EM1} and \cite{EM2} under the assumption that $X$ is metrisable by an ultrametric. Therefore our results are generalisations of some of their results.

We have three applications of Theorem \ref{main theorem}, which connects non-Archimedean analysis and general topology.

First, $\m{C}_{\m{bd}}(X,k)$ satisfies the weak version of the automatic continuity theorem if $k$ is a local field (Theorem \ref{automatic continuity}). Namely, for a Banach $k$-algebra $\A$, every injective $k$-algebra homomorphism $\varphi \colon \m{C}_{\m{bd}}(X,k) \hookrightarrow \A$ with closed image is continuous. In particular, it gives a criterion for the continuity of a faithful linear representation of $\m{C}_{\m{bd}}(X,k)$ on a Banach space.

Second, for an extension $K/k$ of complete valuation fields, the ground field extension $\m{BSC}_K(X) \to \m{BSC}_k(X)$ induced by the inclusion $\m{C}_{\m{bd}}(X,k) \hookrightarrow \m{C}_{\m{bd}}(X,K)$ is a homeomorphism (Proposition \ref{the extension of the scalar}). There is another ground field extension $K \hat{\otimes}_k \m{C}_{\m{bd}}(X,k) \to \m{C}_{\m{bd}}(X,K)$ given by the universality of the complete tensor product $\hat{\otimes}_k$ in the category of Banach $k$-algebras. We will see the difference of those two in Theorem \ref{main theorem 2}.

Finally, we show that the natural continuous map $X \to \m{UF}(X)$ is a homeomorphism onto the image if and only if $X$ is zero-dimensional and Hausdorff (Lemma \ref{completely regular Hausdorff}). We establish Gel'fand theory for totally disconnected compact Hausdorff spaces in this case (Theorem \ref{gelfand}) using a non-Archimedean generalisation of Stone--Weierstrass theorem (\cite{Ber1} 9.2.5.\ Theorem). Here, Gel'fand theory means a natural contravariant-functorial one-to-one correspondence between the collection $\mathscr{C}(X)$ of equivalence classes of totally disconnected compact Hausdorff spaces which contain $X$ as a dense subspace and the set $\mathscr{C}'(X)$ of closed $k$-subalgebras of $\m{C}_{\m{bd}}(X,k)$ separating points of $X$.

\vspace{0.1in}
We remark that the Berkovich spectrum of a Banach algebra is analogous to the Gel'fand transform of a commutative $C^*$-algebra. We study Berkovich spectra in this paper expecting that many facts for Gel'fand transforms also hold for Berkovich spectra. For example, it is well-known that an initial object in the category of compact Hausdorff spaces under $X$ exists and is constructed as the Gel'fand transform $\M_{\C}(\m{C}_{\m{bd}}(X,\C))$ of the commutative $C^*$-algebra $\m{C}_{\m{bd}}(X,\C)$ of bounded continuous $\C$-valued functions on $X$. Therefore our result for the universality of $\m{BSC}_k(X)$ is a direct analogue of this fact. We recall another construction of an initial object in the category of compact Hausdorff spaces under $X$. The Stone--$\check{\m{C}}$ech compactification $\beta X$ of $X$ is constructed as a closed subspace of a direct product of copies of the closed unit disc $\C^{\circ} \subset \C$, and it admits a canonical continuous map $X \to \beta X$ such that every bounded continuous $\C$-valued function on $X$ uniquely extends to a continuous function on $\beta X$. This extension property guarantees that $\beta X$ is also an initial object in the category of compact Hausdorff spaces under $X$. One sometimes assumes that $X$ is a completely regular Hausdorff space in the definition of $\beta X$ so that $X \to \beta X$ is a homeomorphism onto the image, but we do not because we allow compactifications of $X$ whose structure morphism is not injective. Imitating the construction of $\beta X$, we construct a compactification $\m{SC}_k(X)$ of $X$ as a closed subspace of a direct product of copies of the closed unit disc $k^{\circ} \subset k$. We also compare $\m{BSC}_k(X)$ and $\m{SC}_k(X)$, and prove that they are naturally homeomorphic to each other under $X$ when $k$ is a local field or a finite field.

\vspace{0.1in}
In \S \ref{Berkovich Spectra}, we recall the definition of Berkovich spectra. In \S \ref{Stone Spaces}, we recall the Stone space $\m{UF}(X)$ associated to $X$. In \S \ref{Universality of the Stone Space}, we show the universality of $\m{UF}(X)$.

In \S \ref{Statement of the Main Theorem}, we state the main theorem (Theorem \ref{main theorem}). In order to verify it, we construct two set-theoretical maps $\m{supp} \colon \m{BSC}_k(X) \to \m{Spec}(\m{C}_{\m{bd}}(X,k))$ and $\m{Ch}_{\bullet} \colon \m{Spec}(\m{C}_{\m{bd}}(X,k)) \to \m{UF}(X)$. We show that the composite $\m{Ch}_{\m{supp}} \coloneqq \m{Ch}_{\bullet} \circ \m{supp} \colon \m{BSC}_k(X) \to \m{UF}(X)$ is a homeomorphism. Its proof is not straightforward, and is completed in the following two subsections. In \S \ref{Maximality of a Closed Prime Ideal}, we show that every closed prime ideal of $\m{C}_{\m{bd}}(X,k)$ is maximal. In \S \ref{Proof of the Main Theorem}, we verify that the image of $\m{supp}$ coinsides with the subset of closed prime ideals, and we prove that the restriction of $\m{Ch}_{\bullet}$ on the image of $\m{supp}$ is bijective. After that, we verify that $\m{Ch}_{\m{supp}}$ is a homeomorphism, and this completes the proof of Theorem \ref{main theorem}.

In \S \ref{Another Construction}, we compare $\m{BSC}_k(X)$ and $\m{SC}_k(X)$ in the case where $k$ is a local field or a finite field. In \S \ref{Relation to the Stone--Cech compactification}, we observe a connection between $\m{BSC}_{\Q_p}(X)$ and $\beta X$. We show that $\m{BSC}_{\Q_p}(X)$ is homeomorphic to $\beta X$ for special $X$'s.

In \S \ref{Applications}, we deal with the three applications of Theorem \ref{main theorem} mentioned above.

\section{Preliminaries}
\label{Preliminaries}

In this section, we recall the definition of the Berkovich spectrum $\M_k(\A)$ of a Banach algebra $\A$, and the Stone space $\m{UF}(X)$ associated to $X$. For more details, see \cite{Ber1} and \cite{Ber2} for Berkovich spectra, and see  \cite{Ban}, \cite{Joh}, \cite{Sto2}, and \cite{Sto3} for Stone spaces.

\subsection{Berkovich Spectra}
\label{Berkovich Spectra}

A {\it Banach $k$-algebra} means a pair $(\A,\| \cdot \|)$ of a unital associative commutative $k$-algebra $\A$ and a complete submultiplicative non-Archimedean norm $\| \cdot \| \colon \A \to \lbrack 0, \infty )$. We often write $\A$ instead of $(\A,\| \cdot \|)$ for short. Let $(\A,\| \cdot \|)$ be a Banach $k$-algebra. Since $\A$ is unital, it admits a canonical ring homomorphism $k \to \A$, and we also denote by $a \in \A$ the image of $a \in k$. A map $x \colon \A \to \lbrack 0,\infty )$ is said to be a {\it bounded multiplicative seminorm} of $(\A, \| \cdot \|)$ if the following conditions hold:
\begin{itemize}
\item[(i)] $x(f - g) \leq \max \{ x(f), x(g) \}$ for any $f,g \in \A$.
\item[(ii)] $x(fg) = x(f) x(g)$ for any $f,g \in \A$.
\item[(iii)] $x(f) \leq \| f \|$ for any $f \in \A$.
\item[(iv)] $x(a) = |a|$ for any $a \in k$.
\end{itemize}
We denote by $\M_k(\A) = \M_k(\A,\| \cdot \|)$ the set of bounded multiplicative seminorms of $(\A, \| \cdot \|)$ endowed with the weakest topology for which for any $f \in \A$, the map
\begin{eqnarray*}
  f^* \colon \M_k(\A) &\to& \lbrack 0,\infty)\\
  x &\mapsto& x(f)
\end{eqnarray*}
is continuous. We call $\M_k(\A)$ the Berkovich spectrum of $(\A,\| \cdot \|)$. By \cite{Ber1} 1.2.1.\ Theorem, $\M_k(\A)$ is a compact Hausdorff space, and is non-empty if and only if $\A \neq 0$.

\subsection{Stone Spaces}
\label{Stone Spaces}

A $U \subset X$ is said to be {\it clopen} if it is closed and open. We denote by $\m{CO}(X) \subset 2^X$ the set of clopen subsets of $X$. A topological space $X$ is said to be {\it zero-dimensional} if $\m{CO}(X)$ forms an open basis of $X$. The space $\m{CO}(X)$ possesses much information about the topology of $X$ when $X$ is zero-dimensional. The most elementary example of a zero-dimensional space is the underlying topological space of $k$. For each $c \in k$ and $\epsilon > 0$, the subsets of $k$ of the forms $\set{c' \in k}{|c' - c| < \epsilon}$, $\set{c' \in k}{|c' - c| \leq \epsilon}$, $\set{c' \in k}{|c' - c| > \epsilon}$, and $\set{c' \in k}{|c' - c| \geq \epsilon}$ are clopen.

\vspace{0.1in}
The set $\m{CO}(X)$ is a Boolean algebra with respect to $\vee$, $\wedge$, $\neg$, and $\perp$ given by setting $U \vee V \coloneqq U \cup V$, $U \wedge V \coloneqq U \cap V$, $\neg U \coloneqq X \backslash U$, and $\perp \coloneqq \emptyset$ respectively for $U,V \in \m{CO}(X)$. We recall the notion of an ultrafilter of a Boolean algebra. For readers who are not familiar with Boolean algebras and filters, \cite{Joh} and \cite{Sto3} might be helpful. For a Boolean algebra $(A,\vee,\wedge,\neg)$, an $\mathscr{F} \subset A$ is said to be a {\it filter} of $(A,\vee,\wedge,\neg)$ if it satisfies the following:
\begin{itemize}
\item[(i)] $\neg \perp \in \mathscr{F}$.
\item[(ii)] $a \wedge b \in \mathscr{F}$ for any $a,b \in \mathscr{F}$.
\item[(iii)] $a \vee b \in \mathscr{F}$ for any $a \in A$ and $b \in \mathscr{F}$.
\end{itemize}
A filter $\mathscr{F}$ of $(A,\vee,\wedge,\neg)$ is said to be an {\it ultrafilter} if $\mathscr{F} \subsetneq A$ and if for any filter $\mathscr{F}'$ of $(A,\vee,\wedge,\neg)$, $\mathscr{F} \subset \mathscr{F}' \subsetneq A$ implies $\mathscr{F} = \mathscr{F}'$. It is equivalent to the condition that $\perp \notin \mathscr{F}$ and either $a \in \mathscr{F}$ or $\perp a \in \mathscr{F}$ holds for any $a \in A$. For each $S \subset A$, the smallest filter $\mathscr{F}$ of $(A,\vee,\wedge,\neg)$ containing $S$ exists. Then $\mathscr{F}$ is a proper subset of $A$ if and only if $a_1 \wedge \cdots \wedge a_n \neq \perp$ for any $n \in \N \backslash \{ 0 \}$ and $(a_1,\ldots,a_n) \in A^n$. For any filter $\mathscr{F}$ of $(A,\vee,\wedge,\neg)$ with $\mathscr{F} \subsetneq A$, there exists an ultrafilter $\mathscr{F}'$ of $(A,\vee,\wedge,\neg)$ containing $\mathscr{F}$ by Boolean prime ideal theorem. The set of ultrafilters of $(A,\vee,\wedge,\neg)$ is endowed with the topology described in the following way: Its subset $\mathscr{U}$ is open if and only if for any $\mathscr{F} \in \mathscr{U}$, there is an $a \in \mathscr{F}$ such that $\mathscr{G} \in \mathscr{U}$ for any ultrafilter $\mathscr{G}$ of $(A,\vee,\wedge,\neg)$ containing $a$. Applying this construction to $\m{CO}(X)$, we denote by $\m{UF}(X)$ the resulting topological space, and we call it the {\it Stone space} associated to $X$. For example, the subset
\begin{eqnarray*}
  \mathscr{F}(x) \coloneqq \Set{U \in \m{CO}(X)}{x \in U} \subset \m{CO}(X)
\end{eqnarray*}
is an ultrafilter of $(\m{CO}(X),\vee,\wedge,\neg)$ for any $x \in X$, and we call such an ultrafilter a {\it principal ultrafilter}. 

\subsection{Universality of the Stone Space}
\label{Universality of the Stone Space}

We denote by $\m{C}(X,Y)$ the set of continuous maps $f \colon X \to Y$ for topological spaces $X$ and $Y$, and by $\m{Top}$ the category of topological spaces and continuous maps. We also deal with the full subcategory $\m{TDCHTop} \subset \m{Top}$ of totally disconnected compact Hausdorff spaces.

\begin{dfn}
\label{under}
For a category $\mathscr{C}$, a full subcategory $\mathscr{C}' \subset \mathscr{C}$, and an $A \in \m{ob}(\mathscr{C})$, a {\it $\mathscr{C}'$-object under $A$} is a pair $(B,f)$ of a $B \in \m{ob}(\mathscr{C}')$ and an $f \in \m{Hom}_{\mathscr{C}}(A,B)$. Here we regard $B$ as an object of $\mathscr{C}$ through the inclusion $\mathscr{C}' \hookrightarrow \mathscr{C}$. We call $f$ the structure morphism of $(B,f)$ or simply of $B$. We denote by $A/\mathscr{C}'$ the category of $\mathscr{C}'$-objects under $A$ and morphisms compatible with the structure morphisms.
\end{dfn}

In the case $\mathscr{C} = \m{ob}(\m{Top})$, for an $X \in \m{ob}(\m{Top})$ and a $(Y,f) \in \m{ob}(X/\mathscr{C}')$, we call $f$ the structure map of $Y$. We often abbreviate $(Y,f)$ to $Y$.

\begin{dfn}
\label{initial}
For a category $\mathscr{C}$, an object $I$ of $\mathscr{C}$ is said to be {\it initial} if $\m{Hom}_{\mathscr{C}}(I,A)$ consists of one morphism for any $A \in \m{ob}(\mathscr{C})$.
\end{dfn}

An initial object is unique up to a unique isomorphism if it exists. For example, for a category $\mathscr{C}$, a full subcategory $\mathscr{C}' \subset \mathscr{C}$, and an $A \in \m{ob}(\mathscr{C})$, a $(B,\iota) \in \m{ob}(A/\mathscr{C}')$ is initial if and only if the map
\begin{eqnarray*}
  \m{Hom}_{\mathscr{C}}(A,C) & \to & \m{Hom}_{\mathscr{C}'}(B,C) \\
  f & \mapsto & f \circ \iota
\end{eqnarray*}
is bijective for any $C \in \m{ob}(\mathscr{C}')$. In other words, $B$ is an initial $\mathscr{C}'$-object under $A$ with respect to $\iota$ if and only if for any $C \in \m{ob}(\mathscr{C}')$ and any $g \colon A \to C$, there exists a unique $\tilde{g} \colon B \to C$ such that $g = \tilde{g} \circ \iota$.

\begin{thm}
\label{model}
The correspondence $X \rightsquigarrow \m{UF}(X)$ gives a functor $\m{UF} \colon \m{Top} \to \m{TDCHTop}$ which is the left adjoint functor of the inclusion $\m{TDCHTop} \hookrightarrow \m{Top}$.
\end{thm}

We remark that \cite{BJ} Proposition 5.7.12 and the universality of the Stone--$\check{\m{C}}$ech compactification imply Theorem \ref{model}. We will prove Theorem \ref{model} in an explicit way at the end of this subsection. For the proof, we prepare several lemmas and a proposition. We note that for a category $\mathscr{C}$ and a full subcategory $\mathscr{C}' \subset \mathscr{C}$, a functor $J \colon \mathscr{C} \to \mathscr{C}'$ is a left adjoint functor of the inclusion $I \colon \mathscr{C}' \hookrightarrow \mathscr{C}$ if and only if there is a natural transform $\iota \colon \m{id}_{\mathscr{C}} \to I \circ J$ such that the induced map
\begin{eqnarray*}
  \m{Hom}_{\mathscr{C}'}(J(A),B) & \to & \m{Hom}_{\mathscr{C}}(A,I(B)) \\
  f & \mapsto & f \circ \iota_A
\end{eqnarray*}
is bijective for any $A \in \m{ob}(\mathscr{C})$ and $B \in \m{ob}(\mathscr{C}')$. This is equivalent to the condition that $J(A)$ is initial in $A/\mathscr{C}'$ with respect to the adjunction $\iota_A \colon A \to I(J(A)) = J(A)$ for any $A \in \m{ob}(\mathscr{C})$. In order to give a proof of Theorem \ref{model}, we show several fundamental properties of the Stone Space. We remark that this gives an alternative proof of Theorem 3.13 in \cite{Tar}.

\vspace{0.1in}
An $x \in X$ is said to be a {\it cluster point} of an $\mathscr{F} \in \m{UF}(X)$ if $\mathscr{F}$ contains all clopen neighbourhood of $x$. Each $x \in X$ is a cluster point of the principal ultrafilter $\mathscr{F}(x) \in \m{UF}(X)$. Unlike a set-theoretical ultrafilter, the existence of a cluster point gives a strict restriction to an ultrafilter as is shown in the following lemma. An ultrafilter consists of open subsets, and hence carries more information on the topology of $X$ than a set-theoretical ultrafilter does.

\begin{lmm}
\label{principal}
If an $\mathscr{F} \in \m{UF}(X)$ has a cluster point, then $\mathscr{F}$ is a principal ultrafilter.
\end{lmm}

\begin{proof}
Let $x \in X$ be a cluster point of $\mathscr{F}$. Then $\mathscr{F}$ contains the principal ultrafilter $\mathscr{F}(x)$, and hence coincides with $\mathscr{F}(x)$ by the maximality of an ultrafilter.
\end{proof}

For a non-empty family $\mathscr{F}$ of sets, we set $\bigcap \mathscr{F} \coloneqq \bigcap_{U \in \mathscr{F}} U$. We give an explicit description of the set of cluster points of a filter.

\begin{lmm}
\label{cluster point}
The set of cluster points of an $\mathscr{F} \in \m{UF}(X)$ coincides with $\bigcap \mathscr{F}$.
\end{lmm}

\begin{proof}
For a cluster point $x \in X$ of $\mathscr{F}$, one has $x \in \bigcap \mathscr{F}(x) = \bigcap \mathscr{F}$ by Lemma \ref{principal}. For an $x \in \bigcap \mathscr{F}$, assume that there is a $U \in \m{CO}(X)$ such that $x \in U \notin \mathscr{F}$. Then one obtains $X \backslash U \in \mathscr{F}$, and it contradicts the condition $x \in \bigcap \mathscr{F}$. Thus $x$ is a cluster point of $\mathscr{F}$. 
\end{proof}

\begin{lmm}
\label{non-principal ultrafilter}
If $X$ is a discrete infinite set, then $\m{UF}(X)$ contains a non-principal ultrafilter.
\end{lmm}

\begin{proof}
The cardinality of the set of principal ultrafilters is at most $\# X$, while $\# \m{UF}(X)$ coincides with $2^{2^{\# X}}$ in the case where $X$ is a discrete infinite set by \cite{Eng} 3.6.11.\ Theorem.
\end{proof}

\begin{prp}
\label{ultrafilter criterion}
Suppose that $X$ is zero-dimensional.
\begin{itemize}
\item[(i)] $X$ is compact if and only if every ultrafilter has at least one cluster point.
\item[(ii)] $X$ is Hausdorff if and only if every ultrafilter has at most one cluster point.
\item[(iii)] $X$ is a totally disconnected compact Hausdorff space if and only if every ultrafilter has precisely one cluster point.
\end{itemize}
\end{prp}

The assertion is an analogue of the classical result for set-theoretic ultrafilters, and the following proof imitates the proof of it. For the classical result, see \cite{Eng} 1.6.11.\ Proposition and 3.1.24.\ Theorem.

\begin{proof}
When $X$ is zero-dimensional, $X$ is Hausdorff if and only if $X$ is totally disconnected, and therefore the criteria (i) and (ii) immediately imply the criterion (iii).

If $X$ is compact, an ultrafilter has a cluster point because the intersection $\bigcap \mathscr{F}$ is non-empty by the finite-intersection property of a compact space. On the other hand, suppose that every ultrafilter has at least one cluster point. Assume that $X$ is not compact. Since $X$ is zero-dimensional, there is a clopen covering $\mathscr{U}$ of $X$ which has no finite subcovering. The set $\mathscr{V} \coloneqq \set{U \in \m{CO}(X)}{X \backslash U \in \mathscr{U}}$ of complements satisfies $\bigcap \mathscr{V} = \emptyset$ and any finite intersection of clopen subsets in $\mathscr{V}$ is non-empty. Therefore there is an $\mathscr{F} \in \m{UF}(X)$ containing $\mathscr{V}$. One has $\bigcap \mathscr{F} \subset \bigcap \mathscr{V} = \emptyset$, which contradicts the assumption that every ultrafilter has at least one cluster point by Lemma \ref{cluster point}. Thus $X$ is compact.

If $X$ is Hausdorff, then the continuous map $\mathscr{F}(\cdot) \colon X \to \m{UF}(X)$ is injective because $X$ is zero-dimensional. Suppose that every ultrafilter has at most one cluster point. Assume that $X$ is not Hausdorff. There are two distinct points $x,y \in X$ such that any clopen neighbourhoods of $x$ and $y$ have non-empty intersection. In other words, one has $U \cap V \neq \emptyset$ for any $(U,V) \in \mathscr{F}(x) \times \mathscr{F}(y)$. Take a clopen neighbourhood $U \in \mathscr{F}(x)$ of $x$. By the argument above, one has $X \backslash U \notin \mathscr{F}(y)$, and hence $U \in \mathscr{F}(y)$. It implies $\mathscr{F}(x) \subset \mathscr{F}(y)$, and therefore $\mathscr{F}(x) = \mathscr{F}(y)$ by the maximality of an ultrafilter. Both $x$ and $y$ are two distinct cluster points of $\mathscr{F}(x) = \mathscr{F}(y)$, and it contradicts the assumption that every ultrafilter has at most one cluster point. Thus $X$ is Hausdorff.
\end{proof}

As a consequence, for a zero-dimensional space $X$, one obtains the following criteria.
\begin{itemize}
\item[(i)$'$] The space $X$ is compact if and only if $\mathscr{F}(\cdot)$ is surjective.
\item[(ii)$'$] The space $X$ is Hausdorff if and only if $\mathscr{F}(\cdot)$ is injective.
\item[(iii)$'$] The space $X$ is a totally disconnected compact Hausdorff space if and only if $\mathscr{F}(\cdot)$ is bijective.
\end{itemize}
We remark that the bijectivity of $\mathscr{F}(\cdot)$ in (iii)$'$ can be replaced by the condition that $\mathscr{F}(\cdot)$ is a homeomorphism by the following three lemmas.

\begin{lmm}
\label{UF and TDC1}
The $\mathscr{F}(\cdot) \colon X \to \m{UF}(X)$ is continuous and its image is dense.
\end{lmm}

\begin{proof}
For a $U \in \m{CO}(X)$, the pre-image of the open subset $\set{\mathscr{F} \in \m{UF}(X)}{U \in \mathscr{F}}$ is $U \subset X$ itself. Therefore $\mathscr{F}(\cdot)$ is continuous. Let $\mathscr{F} \subset \m{CO}(X)$ be an ultrafilter, and $\mathscr{U} \subset \m{UF}(X)$ an open neighbourhood of $\mathscr{F}$. By the definition of the topology of $\m{UF}(X)$, there is a $U \in \m{CO}(X)$ such that $U \in \mathscr{F}$ and $\mathscr{V} \coloneqq \set{\mathscr{F}' \in \m{UF}(X)}{U \in \mathscr{F}'} \subset \mathscr{U}$. Then $U \neq \emptyset$ because $\emptyset \notin \mathscr{F}$, and hence $\mathscr{F}(U) \neq \emptyset$. Since $\mathscr{F}(U) \subset \mathscr{F}(X) \cap \mathscr{V} ap \subset \mathscr{F}(X) \cap \mathscr{U}$, one concludes $\mathscr{F}(X) \cap \mathscr{U} \neq \emptyset$.
\end{proof}

\begin{lmm}
\label{UF and TDC2}
The space $\m{UF}(X)$ is a totally disconnected compact Hausdorff space.
\end{lmm}

This assertion is contained in the general fact of the Stone space in \cite{Sto2} Theorem IV${}_2$, but we give a proof for reader's convenience.

\begin{proof}
For a $U \in \m{CO}(X)$, one has
\begin{eqnarray*}
  \m{UF}(X) = \Set{\mathscr{F} \in \m{UF}(X)}{U \in \mathscr{F}} \sqcup \Set{\mathscr{F} \in \m{UF}(X)}{X \backslash U \in \mathscr{F}},
\end{eqnarray*}
and hence $\m{CO}(\m{UF}(X))$ forms an open basis of $\m{UF}(X)$. Therefore by Proposition \ref{ultrafilter criterion} and Lemma \ref{UF and TDC1}, it suffices to show that $\m{UF}(X)$ is compact and Hausdorff, because a continuous map from a compact space to a Hausdorff space is a closed map.

\vspace{0.1in}
For $\mathscr{F}, \mathscr{G} \in \m{UF}(X)$ with $\mathscr{F} \neq \mathscr{G}$, take a $U \in \m{CO}(X)$ contained in precisely one of them. Then the complement $X \backslash U$ is contained in the other one. Therefore the partition
\begin{eqnarray*}
  \m{UF}(X) = \Set{\mathscr{F} \in \m{UF}(X)}{U \in \mathscr{F}} \sqcup \Set{\mathscr{F} \in \m{UF}(X)}{X \backslash U \in \mathscr{F}}
\end{eqnarray*}
by clopen subsets of $\m{UF}(X)$ separates $\mathscr{F}$ and $\mathscr{G}$. Thus $\m{UF}(X)$ is Hausdorff.

\vspace{0.1in}
Assume that $\m{UF}(X)$ is not compact. There is a clopen covering $\mathscr{U}$ of $\m{UF}(X)$ which has no finite subcovering. In particular, the subset
\begin{eqnarray*}
  \mathscr{V} \coloneqq \Set{U \in \m{CO}(\m{UF}(X))}{\m{UF}(X) \backslash U \in \mathscr{U}}
\end{eqnarray*}
satisfies $\bigcap \mathscr{V} = \emptyset$ and any finite intersection of clopen subsets belonging to $\mathscr{V}$ is non-empty. Since the map $\mathscr{F}(\cdot)$ is continuous, the inverse image
\begin{eqnarray*}
  \mathscr{F}(\cdot)^* \mathscr{V} \coloneqq \Set{\mathscr{F}(\cdot)^{-1}(V)}{V \in \mathscr{V}}
\end{eqnarray*}
is a non-empty subset of $\m{CO}(X)$ satisfying that $\bigcap \mathscr{F}(\cdot)^* \mathscr{V} = \emptyset$ and any finite intersection of clopen subsets belonging to $\mathscr{F}(\cdot)^* \mathscr{V}$ is non-empty. Therefore there is an $\mathscr{F} \in \m{UF}(X)$ containing $\mathscr{F}(\cdot)^* \mathscr{V}$ by the facts recalled in \S \ref{Stone Spaces}. Since $\mathscr{U}$ covers $\m{UF}(X)$, there is a $U \in \mathscr{U}$ containing $\mathscr{F}$. The pre-image $V \in \mathscr{F}(\cdot)^* \mathscr{V}$ of the complement $\m{UF}(X) \backslash U \in \mathscr{V}$ is contained in $\mathscr{F}$ because $\mathscr{F}(\cdot)^* \mathscr{V} \subset \mathscr{F}$. By the definition of the topology of $\m{UF}(X)$, there is a $W \in \mathscr{F}$ such that $W \in \mathscr{G}$ implies $\mathscr{G} \in U$ for any $\mathscr{G} \in \m{UF}(X)$. In particular, for any $x \in W$, $W \subset \mathscr{F}(x)$ and hence $\mathscr{F}(x) \in U$. Therefore one obtains $W \subset \mathscr{F}(\cdot)^{-1}(U)$. Since $V,W \in \mathscr{F}$, one has $V \cap W \in \mathscr{F}$ and hence $V \cap W \neq \emptyset$. Take an $x \in V \cap W \subset X$. Since $V = \mathscr{F}(\cdot)^{-1}(\m{UF}(X) \backslash U)$, one has $\mathscr{F}(x) \notin U$, which contradicts the condition $x \in W \subset \mathscr{F}(\cdot)^{-1}(U)$. Thus  $\m{UF}(X)$ is compact.
\end{proof}

\begin{lmm}
\label{reflexivity}
If $X$ is a totally disconnect compact Hausdorff space, then $\mathscr{F}(\cdot) \colon X \to \m{UF}(X)$ is a homeomorphism.
\end{lmm}

In particular, $\mathscr{F}(\cdot) \colon \m{UF}(X) \to \m{UF}(\m{UF}(X))$ is a homeomorphism without the assumption on $X$ by Lemma \ref{UF and TDC2}.

\begin{proof}
The assertion immediately follows from Proposition \ref{ultrafilter criterion} (iii), Lemma \ref{UF and TDC1}, and Lemma \ref{UF and TDC2}, because every continuous map between compact Hausdorff spaces is closed.
\end{proof}

\begin{proof}[Proof of Theorem \ref{model}]
By Lemma \ref{UF and TDC1} and Lemma \ref{UF and TDC2}, $(\m{UF}(X),\mathscr{F}(\cdot))$ is an object of $X/\m{TDCHTop}$. Let $Y,Z \in \m{Top}$ and $f \in \m{C}(Y,Z)$. For an $\mathscr{F} \in \m{UF}(Y)$, the subset
\begin{eqnarray*}
  \m{UF}(f)_* \mathscr{F} \coloneqq \Set{U \in \m{CO}(Y)}{\varphi^{-1}(U) \in \mathscr{F}}.
\end{eqnarray*}
is an ultrafilter of $\m{CO}(Z)$. The map $\m{UF}(f)_* \colon \m{UF}(Y) \to \m{UF}(Z)$ is continuous by the definition of the topologies of $\m{UF}(Y)$ and $\m{UF}(Z)$. The correspondences $Y \rightsquigarrow \m{UF}(Y)$ and $f \rightsquigarrow \m{UF}(f)_*$ gives a functor $\m{UF} \colon \m{Top} \to \m{TDCHTop}$. Therefore it suffices to show that $(\m{UF}(X),\mathscr{F}(\cdot))$ is an initial object of $X/\m{TDCHTop}$.

Let $(Y,\varphi)$ be an object of $X/\m{TDCHTop}$. Since the image of $X$ is dense in $\m{UF}(X)$ by Lemma \ref{UF and TDC1} and $Y$ is Hausdorff, a continuous extension $\m{UF}(\varphi) \colon \m{UF}(X) \to Y$ is unique if it exists. The diagram
\begin{eqnarray*}
\begin{CD}
  X         @>{\varphi}>>           Y \\
  @V{\mathscr{F}(\cdot)}VV          @VV{\mathscr{F}(\cdot)}V \\
  \m{UF}(X) @>{\m{UF}(\varphi)_*}>> \m{UF}(Y)
\end{CD}
\end{eqnarray*}
commutes by the definitions of $\mathscr{F}(\cdot)$ and $\m{UF}(\varphi)_*$, and the right vertical map is a homeomorphism by Lemma \ref{reflexivity}. Therefore one obtains a continuous extension $\mathscr{F}(\cdot)^{-1} \circ \m{UF}(\varphi)_* \colon \m{UF}(X) \to Y$ of $\varphi$.
\end{proof}

\section{Main Result}
\label{Main Result}

\subsection{Statement of the Main Theorem}
\label{Statement of the Main Theorem}

We denote by $\m{C}_{\m{bd}}(X,k)$ the Banach $k$-algebra of bounded continuous $k$-valued functions on $X$ endowed with the supremum norm. We put $\m{BSC}_k(X) \coloneqq \M_k(\m{C}_{\m{bd}}(X,k))$. Let $\iota_k$ denote the evaluation map
\begin{eqnarray*}
  \iota_k \colon X & \to & \m{BSC}_k(X) \\
  x & \mapsto & \left( \iota_k(x) \colon f \mapsto |f(x)| \right),
\end{eqnarray*}
which is continuous by the definition of the topology of $\m{BSC}_k(X)$.

\begin{thm}
\label{main theorem}
There is a natural homeomorphism $\m{BSC}_k(X) \cong \m{UF}(X)$ compatible with $\iota_k$ and $\mathscr{F}(\cdot)$.
\end{thm}

In other words, there is a natural transform $\Phi \colon \m{BSC}_k \to \m{UF}$ such that $\Phi(Y)$ lies in $\m{Hom}_{Y/\m{TDCHTop}}((\m{BSC}(Y),\iota_k),(\m{UF}(Y),\mathscr{F}(\cdot)))$ for any topological space $Y$. In particular, it gives an isomorphism $(\m{BSC}(X),\iota_k) \cong (\m{UF}(X),\mathscr{F}(\cdot))$ in $X/\m{TDCHTop}$, and hence $(\m{BSC}(X),\iota_k)$ satisfies the same universality as $(\m{UF}(X),\mathscr{F}(\cdot))$ does.

\begin{crl}
\label{initial 2}
The space $\m{BSC}_k(X)$ is initial in $X/\m{TDCHTop}$ with respect to $\iota_k$.
\end{crl}

\begin{crl}
The functor
\begin{eqnarray*}
  \m{BSC}_k \colon \m{Top} & \to & \m{TDCHTop} \\
  X & \rightsquigarrow & \m{BSC}_k(X)
\end{eqnarray*}
is a left adjoint functor of the inclusion of the full subcategory.
\end{crl}

\begin{crl}
\label{dense}
The image of $\iota_k \colon X \to \m{BSC}_k(X)$ is dense.
\end{crl}

In order to prove Theorem \ref{main theorem}, we introduce two set-theoretical maps $\m{supp}$ and $\m{Ch}_{\bullet}$. For an $x \in \m{BSC}_k(X)$, its support $\m{supp}(x) \coloneqq \set{f \in \m{C}_{\m{bd}}(X,k)}{x(f) = 0}$ is a closed prime ideal. We call the map
\begin{eqnarray*}
  \m{supp} \colon \m{BSC}_k(X) &\to& \m{Spec}(\m{C}_{\m{bd}}(X,k)) \\
  x & \mapsto & \m{supp}(x)
\end{eqnarray*}
the support map. For an $m \in \m{Spec}(\m{C}_{\m{bd}}(X,k))$, the family $\m{Ch}_m \coloneqq \set{U \in \m{CO}(X)}{1_U \notin m}$ is an ultrafilter, where $1_U \colon X \to k$ denotes the characteristic function of $U \in \m{CO}(X)$. Indeed, $\m{Ch}_m$ is stable under $\cup$ because $m$ is an ideal, and is stable under $\cap$ because $m$ is a prime ideal. The maximality of $\m{Ch}_m$ follows from the property that either $1_U \in m$ or $1_{X \backslash U} = 1 - 1_U \in m$ holds for any $U \in \m{CO}(X)$ because $m$ is a prime ideal. We call the map
\begin{eqnarray*}
  \m{Ch}_{\bullet} \colon \m{Spec}(\m{C}_{\m{bd}}(X,k)) & \to & \m{UF}(X) \\
  m & \mapsto & \m{Ch}_m
\end{eqnarray*}
the characteristic map. We put $\m{Ch}_{\m{supp}} \coloneqq \m{Ch}_{\bullet} \circ \m{supp} \colon \m{BSC}_k(X) \to \m{UF}(X)$.

\begin{exm}
\label{compatibility}
For an $x \in X$, $\m{supp}(\iota_k(x)) \subset \m{C}_{\m{bd}}(X,k)$ is the maximal ideal consisting of functions vanishing at $x$, and one has $\m{Ch}_{\m{supp}(\iota_k(x))} = \mathscr{F}(x)$. Thus $\m{Ch}_{\m{supp}}$ is an extension of the continuous map $\mathscr{F}(\cdot) \colon X \to \m{UF}(X)$ via $\iota_k$.
\end{exm}

We prove that $\m{Ch}_{\m{supp}}$ is a homeomorphism under $X$ in three steps in \S \ref{Maximality of a Closed Prime Ideal} and \S \ref{Proof of the Main Theorem}. First, we show that every closed prime ideal of $\m{C}_{\m{bd}}(X,k)$ is a maximal ideal. Second, we verify that the image of $\m{supp}$ coincides with the subset of closed prime ideals, and study the restriction of $\m{Ch}_{\bullet}$ on the image of $\m{supp}$. Finally, we prove that $\m{Ch}_{\m{supp}}$ is a homeomorphism.

\subsection{Maximality of a Closed Prime Ideal}
\label{Maximality of a Closed Prime Ideal}

We prove that every closed prime ideal of $\m{C}_{\m{bd}}(X,k)$ is a maximal ideal. We remark that this is proved by Alain Escassut and Nicolas Ma\"inetti in \cite{EM1} Theorem 12 in the case where $X$ is an ultrametric space. Here we assume nothing on $X$, and hence $X$ is not necessarily metrisable.

\begin{prp}
\label{anti-order preserving}
For any $m_1, m_2 \in \m{Spec}(\m{C}_{\m{bd}}(X,k))$ with $m_1 \subset m_2$, $\m{Ch}_{m_1} = \m{Ch}_{m_2}$.
\end{prp}

\begin{proof}
The condition $m_1 \subset m_2$ implies $\m{Ch}_{m_2} \subset \m{Ch}_{m_1}$ by definition. Since $\m{Ch}_{m_2}$ is an ultrafilter, the inclusion guarantees $\m{Ch}_{m_1} = \m{Ch}_{m_2}$.
\end{proof}

\begin{prp}
\label{closed prime and ultrafilter}
For closed prime ideals $m_1, m_2 \subset \m{C}_{\m{bd}}(X,k)$, the equality $\m{Ch}_{m_1} = \m{Ch}_{m_2}$ implies $m_1 = m_2$.
\end{prp}

\begin{proof}
Suppose $\m{Ch}_{m_1} = \m{Ch}_{m_2}$ for closed prime ideals $m_1, m_2 \subset \m{C}_{\m{bd}}(X,k)$. It suffices to show $m_1 \subset m_2$. Take an element $f \in m_1$. For a positive real number $\epsilon$, we set $U_{\epsilon} \coloneqq \set{x \in X}{|f(x)| < \epsilon}$, and then $U_{\epsilon} \subset X$ is a clopen subset, because it is preimage of the clopen subset $\set{c \in k}{|f(x) - c| < \epsilon}$ by the continuous function $f$. Set $f_{\epsilon} \coloneqq (1 - 1_{U_{\epsilon}})f \in \m{C}_{\m{bd}}(X,k)$. Since $f \in m_1$, one has $f_{\epsilon} \in m_1$. The absolute value of $f_{\epsilon} + 1_{U_{\epsilon}} \in \m{C}_{\m{bd}}(X,k)$ at each point in $X$ has a lower bound $\min \{ \epsilon, 1 \}$, and hence its inverse is bounded and continuous. It implies that $f_{\epsilon} + 1_{U_{\epsilon}}$ is invertible in $\m{C}_{\m{bd}}(X,k)$, and therefore $1_{U_{\epsilon}} \notin m_1$. One has $U_{\epsilon} \in \m{Ch}_{m_1} = \m{Ch}_{m_2}$, and hence $1 - 1_{U_{\epsilon}} = 1_{X \backslash U_{\epsilon}} \in m_2$. Thus $f_{\epsilon} = (1 - 1_{U_{\epsilon}})f \in m_2$, and the inequality $\| f - f_{\epsilon} \| = \| 1_{U_{\epsilon}}f \| \leq \epsilon$ guarantees $f \in m_2$ by the closedness of $m_2$.
\end{proof}

\begin{prp}
\label{closed prime and maximal}
Every closed prime ideal of $\m{C}_{\m{bd}}(X,k)$ is a maximal ideal.
\end{prp}

We note that for a Banach $k$-algebra $\A$, every maximal ideal of $\A$ is a closed prime ideal by \cite{BGR} 1.2.4.\ Corollary 5, but the converse does not hold in general. For example, the Tate algebra $k \{ T \}$ has a non-maximal closed ideal $\{ 0 \} \subset k \{ T \}$.

\begin{proof}
For a closed prime ideal $m_1 \subset \m{C}_{\m{bd}}(X,k)$, take a maximal ideal $m_2 \subset \m{C}_{\m{bd}}(X,k)$ containing $m_1$. Then $m_2$ is also a closed prime ideal by \cite{BGR} 1.2.4.\ Corollary 5. The assertion immediately follows from Proposition \ref{anti-order preserving} and Proposition \ref{closed prime and ultrafilter}.
\end{proof}

\subsection{Proof of the Main Theorem}
\label{Proof of the Main Theorem}

\begin{prp}
\label{image of supp}
The image of $\m{supp}$ is the subset of closed prime ideals.
\end{prp}

\begin{proof}
Every closed prime ideal $m \subset \m{C}_{\m{bd}}(X,k)$ is a maximal ideal by Proposition \ref{closed prime and maximal}, and hence there is an $x \in \m{BSC}_k(X)$ such that $\m{supp}(x) = m$ by the argument in the proof of \cite{Ber1} 1.2.1.\ Theorem.
\end{proof}

\begin{prp}
\label{maximal and utrafilter}
The restriction of $\m{Ch}_{\bullet}$ on the image of $\m{supp}$ is bijective.
\end{prp}

\begin{proof}
If $X = \emptyset$, then $\m{Spec}(\m{C}_{\m{bd}}(X,k)) = \m{UF}(X) = \emptyset$, and hence we may assume $X \neq \emptyset$. By Proposition \ref{closed prime and ultrafilter} and Proposition \ref{image of supp}, it suffices to verify the surjectivity. Take an $\mathscr{F} \in \m{UF}(X)$. Set
\begin{eqnarray*}
  m \coloneqq \Set{f \in \m{C}_{\m{bd}}(X,k)}{\inf_{U \in \mathscr{F}} \sup_{x \in U} |f(x)| = 0} \subset \m{C}_{\m{bd}}(X,k).
\end{eqnarray*}
Then $m \subset \m{C}_{\m{bd}}(X,k)$ is an ideal, and $1 \notin m$ because $|1(x)| = 1$ for any $x \in X \neq \emptyset$. We verify that the map
\begin{eqnarray*}
  \| \cdot \|_{\mathscr{F}} \colon \m{C}_{\m{bd}}(X,k) &\to& \lbrack 0,\infty )\\
  f &\mapsto& \inf_{U \in \mathscr{F}} \sup_{x \in U} |f(x)| < \| f \|
\end{eqnarray*}
is continuous. The map $\| \cdot \|_{\mathscr{F}}$ is continuous at any $f \in \m{C}_{\m{bd}}(X,k)$ with $\| f \|_{\mathscr{F}} = 0$ because for any $g \in \m{C}_{\m{bd}}(X,k) \backslash \{ f \}$, there is a $U_0 \in \mathscr{F}$ with $\sup_{x \in U} |f(x)| < \| f - g \|$ and hence
\begin{eqnarray*}
  \| g \|_{\mathscr{F}} \leq \inf_{U \in \mathscr{F}} \sup_{x \in U} |f(x) - (f - g)(x)| \leq \inf_{U \in \mathscr{F}} \sup_{x \in U} \max \left\{ |f(x)|, |(f - g)(x)| \right\} \leq \| f - g \|.
\end{eqnarray*}
The map $\| \cdot \|_{\mathscr{F}}$ is locally constant at any $f \in \m{C}_{\m{bd}}(X,k)$ with $\| f \|_{\mathscr{F}} \neq 0$ because for any $g \in \m{C}_{\m{bd}}(X,k)$ with $\| f - g \| < \| f \|_{\mathscr{F}}$, we have
\begin{eqnarray*}
  \| g \|_{\mathscr{F}} & \leq & \inf_{U \in \mathscr{F}} \sup_{x \in U} |f(x) - (f - g)(x)| \leq \inf_{U \in \mathscr{F}} \sup_{x \in U} \max \left\{ |f(x)|, |(f - g)(x)| \right\} \\
  & \leq & \inf_{U \in \mathscr{F}} \max \left\{ \sup_{x \in U} |f(x)|, \| f - g \| \right\} = \| f \|_{\mathscr{F}}.
\end{eqnarray*}
Therefore $\| \cdot \|_{\mathscr{F}}$ is continuous. Since $\{ 0 \} \subset \lbrack 0,\infty )$ is closed, $m$ is a closed ideal. For $f,g \in \m{C}_{\m{bd}}(X,k)$ with $fg \in m$, suppose $f \notin m$. We prove $g \in m$. If $g = 0$, then $g \in m$. Therefore we may assume $g \neq 0$. Since $f \notin m$, there is some $\epsilon_0 > 0$ such that the clopen subset $V \coloneqq \set{x \in X}{|f(x)| < \epsilon}$ does not belong to $\mathscr{F}$ for any $0 < \epsilon < \epsilon_0$. Let $0 < \epsilon < \epsilon_0$. The condition $fg \in m$ implies that there is some $U \in \mathscr{F}$ such that $\sup_{x \in U} |(fg)(x)| < \epsilon^2$. Since $\mathscr{F}$ is an ultrafilter, one has $X \backslash V \in \mathscr{F}$ and hence $U \backslash V = U \cap (X \backslash V) \in \mathscr{F}$. For an $x \in U \backslash V$, the inequality $|g(x)| = |f(x)|^{-1} |f(x)g(x)| < \epsilon$ implies $\sup_{x \in U \backslash V} |g(x)| < \epsilon$. One obtains $\| g \|_{\mathscr{F}} = 0$, and hence $g \in m$. Therefore $m$ is a closed prime ideal. Let $U \in \mathscr{F}$. One gets $\| 1_U \|_{\mathscr{F}} = 1$ by definition, and hence $U \in \m{Ch}_m$. It implies $\mathscr{F} \subset \m{Ch}_m$. Since $\mathscr{F}$ is an ultrafilter, one concludes $\mathscr{F} = \m{Ch}_m$. Thus $\mathscr{F}$ is contained in the image of $\m{Ch}_{\m{supp}}$ by Proposition \ref{image of supp}.
\end{proof}

\begin{proof}[Proof of Theorem \ref{main theorem}]
The map $\m{Ch}_{\m{supp}}$ is compatible with $\iota_k$ and $\mathscr{F}(\cdot)$ as is shown in Example \ref{compatibility}. We prove that $\m{Ch}_{\m{supp}}$ is a homeomorphism. We first prove the bijectivity. Since the restriction of $\m{Ch}_{\bullet}$ on the image of $\m{supp}$ is bijective by Proposition \ref{maximal and utrafilter}, we have only to show that $\m{supp}$ is injective. For that purpose, for a maximal ideal $m \subset \m{C}_{\m{bd}}(X,k)$, we consider the relation between the quotient seminorm $\| \cdot + m \|$ at $m$ and the map $\| \cdot \|_{\m{Ch}_m}$ defined in the proof of Proposition \ref{maximal and utrafilter}. For an $f \in \m{C}_{\m{bd}}(X,k)$, one has
\begin{eqnarray*}
  \| f + m \| = \inf_{g \in m} \| f - g \| \geq \inf_{g \in m} \| f - g \|_{\m{Ch}_m} = \inf_{g \in m} \| f \|_{\m{Ch}_m} = \| f \|_{\m{Ch}_m}.
\end{eqnarray*}
Take an $r \in \R$ with $\| f \|_{\m{Ch}_m} < r$. Set
\begin{eqnarray*}
  U \coloneqq \Set{x \in X}{|f(x)| > r}.
\end{eqnarray*}
Then $U \subset X$ is clopen by an argument similar to the one in the proof of Proposition \ref{closed prime and ultrafilter}. If $U \in \m{Ch}_m$, then one has
\begin{eqnarray*}
  \| f \|_{\m{Ch}_m} = \inf_{V \in \m{Ch}_m} \sup_{x \in V} |f(x)| \geq \inf_{V \in \m{Ch}_m} \sup_{x \in U \cap V} |f(x)| \geq \inf_{V \in \m{Ch}_m} r = r
\end{eqnarray*}
and hence it contradicts the condition $\| f \|_{\m{Ch}_m} < r$. It implies $U \notin \m{Ch}_m$, and therefore $1_U \in m$. One obtains
\begin{eqnarray*}
  \| f + m \| \leq \| f - 1_Uf \| = \| 1_{X \backslash U}f \| \leq r.
\end{eqnarray*}
One gets $\| f + m \| = \| f \|_{\m{Ch}_m}$.

\vspace{0.1in}
Next, we prove that the map $\| \cdot \|_{\m{Ch}_m}$ is a bounded multiplicative seminorm on $\m{C}_{\m{bd}}(X,k)$. It is a bounded power-multiplicative seminorm by definition, and it suffices to show the multiplicativity. Let $f,g \in \m{C}_{\m{bd}}(X,k)$ such that $\| fg \|_{\m{Ch}_m} < \| f \|_{\m{Ch}_m} \| g \|_{\m{Ch}_m}$. In particular, $\| f \|_{\m{Ch}_m} \| g \|_{\m{Ch}_m} \neq 0$ and $f,g \notin m$. Take an $\epsilon > 0$ such that $\epsilon < \| f \|_{\m{Ch}_m}$, $\epsilon < \| g \|_{\m{Ch}_m}$, and $\| fg \|_{\m{Ch}_m} < (\| f \|_{\m{Ch}_m} - \epsilon)(\| g \|_{\m{Ch}_m} - \epsilon)$. Set
\begin{eqnarray*}
  && V_1 \coloneqq \Set{x \in X}{|f(x)| > \| f \|_{\m{Ch}_m}  - \epsilon}\\
  && V_2 \coloneqq \Set{x \in X}{|g(x)| > \| g \|_{\m{Ch}_m} - \epsilon}.
\end{eqnarray*}
Then $V_1, V_2 \subset X$ are clopen. If $V_1 \notin \m{Ch}_m$, then $X \backslash V_1 \in \m{Ch}_m$, but the inequality
\begin{eqnarray*}
  \| f \|_{\m{Ch}_m} \leq \sup_{x \in X \backslash V_1} |f(x)| \leq \| f \|_{\m{Ch}_m} - \epsilon
\end{eqnarray*}
contradicts the condition $\epsilon > 0$. Therefore $V_1 \in \m{Ch}_m$. Similarly, one obtains $V_2 \in \m{Ch}_m$, and hence $V_1 \cap V_2 \in \m{Ch}_m$. Then the inequality
\begin{eqnarray*}
  && \| fg \|_{\m{Ch}_m} < (\| f \|_{\m{Ch}_m}  - \epsilon)(\| g \|_{\m{Ch}_m}  - \epsilon) \leq \inf_{W \in \m{Ch}_m} \sup_{x \in V_1 \cap V_2 \cap W} |f(x)| \ |g(x)| \\
  &\leq& \inf_{W \in \m{Ch}_m} \sup_{x \in W} |f(x)g(x)| = \| fg \|_{\m{Ch}_m}
\end{eqnarray*}
holds, and it is a contradiction. Thus $\| fg \|_{\m{Ch}_m} = \| f \|_{\m{Ch}_m} \| g \|_{\m{Ch}_m}$. We conclude that the map $\| \cdot \|_{\m{Ch}_m}$ is a bounded multiplicative seminorm, and hence corresponds to a point in $\m{BSC}_k(X)$.

\vspace{0.1in}
Now take an $x \in \m{BSC}_k(X)$. Since $y \coloneqq \| \cdot \|_{\m{Ch}_{\m{supp}}(x)} \in \m{BSC}_k(X)$ coincides with the quotient seminorm $\| \cdot + \m{supp}(x) \|$, one has $x(f) \leq y(f)$ for any $f \in \m{C}_{\m{bd}}(X,k)$. It implies that $x$ gives a bounded multiplicative norm of the complete residue field $k(y)$ at $y$, because $\m{supp}(y)$ is a maximal ideal. It implies $x = y$ because $y(f) = y(f^{-1})^{-1} \leq x(f^{-1})^{-1} = x(f)$ for any $f \in k(y)^{\times}$. Thus $x$ is reconstructed from its image $y$ by $\m{Ch}_{\m{supp}}$, and hence $\m{Ch}_{\m{supp}}$ is injective.

\vspace{0.1in}
Finally, we verify the continuity of $\m{Ch}_{\m{supp}}$. Take a $U \in \m{CO}(X)$, and set $\mathscr{U} \coloneqq \set{\mathscr{F} \in \m{UF}(X)}{U \in \mathscr{F}}$. The pre-image of $\mathscr{U}$ by $\m{Ch}_{\m{supp}}$ is the subset
\begin{eqnarray*}
  && \Set{x \in \m{BSC}_k(X)}{U \in \m{Ch}_{\m{supp}}(x)} = \Set{x \in \m{BSC}_k(X)}{1_U \notin \m{supp}(x)}\\
  &=& \Set{x \in \m{BSC}_k(X)}{x(1_U) > 0} \subset \m{BSC}_k(X),
\end{eqnarray*}
and it is open by the definition of the topology of $\m{BSC}_k(X)$. Therefore $\m{Ch}_{\m{supp}}$ is a continuous bijective map between compact Hausdorff spaces, and is a homeomorphism. This completes the proof.
\end{proof}

We give several corollaries. These are generalisations of some of results in \cite{EM1} and \cite{EM2}. In those papers, Alain Escassut and Nicolas Ma\"inetti deal with ultrametric spaces, while we deal with general topological spaces. We remark that they deal with not only the class of bounded continuous functions, but also that of bounded uniformly continuous functions with respect to the uniform structure associated to the ultrametric.

\begin{crl}
The map $\m{supp}$ gives a bijective map from $\m{BSC}_k(X)$ to the set of maximal ideals of $\m{C}_{\m{bd}}(X,k)$, and every maximal ideal of $\m{C}_{\m{bd}}(X,k)$ is the support of a unique bounded multiplicative seminorm on $\m{C}_{\m{bd}}(X,k)$.
\end{crl}

This is a generalisation of \cite{EM1} Theorem 16 for the class of bounded continuous functions.

\begin{proof}
We proved that the injectivity of $\m{supp}$ in the proof of Theorem \ref{main theorem}, and the image of $\m{supp}$ coincides with the subset of maximal ideals by Proposition \ref{closed prime and maximal} and Proposition \ref{image of supp}. Thus the assertion holds.
\end{proof}

\begin{crl}
Every bounded multiplicative seminorm on $\m{C}_{\m{bd}}(X,k)$ is of the form
\begin{eqnarray*}
  \| \cdot \|_{\mathscr{F}} \colon \m{C}_{\m{bd}}(X,k) &\to& \lbrack 0,\infty )\\
  f &\mapsto& \inf_{U \in \mathscr{F}} \sup_{x \in U} |f(x)|
\end{eqnarray*}
for a unique $\mathscr{F} \in \m{CO}(X)$.
\end{crl}

\begin{proof}
Let $x \in \m{BSC}_k(X)$. We proved the equality $x = \| \cdot \|_{\m{Ch}_{\m{supp}}(x)}$ in the proof of Theorem \ref{main theorem}. The uniqueness of an $\mathscr{F} \in \m{CO}(X)$ follows from the surjectivity of $\m{Ch}_{\m{supp}}$.
\end{proof}

We denote by $\m{UF}(|X|)$ the set of set-theoretical ultrafilters of $X$. We compare $\m{UF}(|X|)$ with $\m{UF}(X)$ through the bijection $\m{Ch}_{\m{supp}}$ in Theorem \ref{main theorem}.

\begin{crl}
\label{uf - suf}
The inclusion $\m{CO}(X) \hookrightarrow 2^X$ is a Boolean algebra homomorphism, and induces a surjective map
\begin{eqnarray*}
  (\cdot \cap \m{CO}(X)) \colon \m{UF}(|X|) &\to& \m{UF}(X) \\
  \mathscr{U} &\mapsto& \mathscr{U} \cap \m{CO}(X).
\end{eqnarray*}
For $\mathscr{U}, \mathscr{U}' \in \m{UF}(|X|)$, the equality $\lim_{\mathscr{U}} |f(x)| = \lim_{\mathscr{U}'} |f(x)|$ holds for any $f \in \m{C}_{\m{bd}}(X,k)$ if and only if $\mathscr{U} \cap \m{CO}(X) = \mathscr{U}' \cap \m{CO}(X)$.
\end{crl}

\begin{proof}
Let $\mathscr{F} \in \m{UF}(X)$. Since $\mathscr{F}$ is a family of subsets of $X$ which is closed under intersections and satisfies $\emptyset \notin \mathscr{F}$, there is an $\mathscr{U} \in \m{UF}(|X|)$ containing $\mathscr{F}$. It implies the surjectivity of the given correspondence. Let $\mathscr{U} \in \m{UF}(|X|)$ and $f \in \m{C}_{\m{bd}}(X,k)$. The limit $\lim_{\mathscr{U}} |f(x)|$ exists because the boundedness of $f$ guarantees that $f(X)$ is relatively compact in $\R$. Moreover, since $\mathscr{U} \cap \m{CO}(X) \subset \mathscr{U}$, we have $\| f \|_{\mathscr{U} \cap \m{CO}(X)} = \lim_{\mathscr{U}} |f(x)|$. Thus the second assertion follows from the injectivity of the inverse map of $\m{Ch}_{\m{supp}} \colon \m{BSC}_k(X) \to \m{UF}(X)$.
\end{proof}

\begin{crl}
Every bounded multiplicative seminorm on $\m{C}_{\m{bd}}(X,k)$ is of the form
\begin{eqnarray*}
  \m{C}_{\m{bd}}(X,k) &\to& \lbrack 0,\infty )\\
  f &\mapsto& \lim_{\mathscr{U}} |f(x)|
\end{eqnarray*}
for a $\mathscr{U} \in \m{UF}(|X|)$, where $\lim_{\mathscr{U}} |f(x)|$ denotes the limit of the $\R$-valued continuous function $|f| \colon X \to \R \colon x \mapsto |f(x)|$ along $\mathscr{U}$ for each $f \in \m{C}_{\m{bd}}(X,k)$.
\end{crl}

This together with Corollary \ref{uf - suf} is a generalisation of \cite{EM1} Corollary 16.3.

\begin{proof}
Every $x \in \m{BSC}_k(X)$ is presented as $\| \cdot \|_{\mathscr{F}}$ by $\mathscr{F} \coloneqq \m{Ch}_{\m{supp}}(x)$. By Corollary \ref{uf - suf}, there is a $\mathscr{U} \in \m{UF}(|X|)$ containing $\mathscr{F}$, and satisfies $x(f) = \| f \|_{\mathscr{F}} = \lim_{\mathscr{U}} |f(x)|$ for any $f \in \m{C}_{\m{bd}}(X,k)$.
\end{proof}

A topological space $X$ is said to be {\it strongly zero-dimensional} if for any disjoint closed subsets $F, F' \subset X$ there is a $U \in \m{CO}(X)$ such that $F \subset U \subset X \backslash F'$. We note that every strongly zero-dimensional Hausdorff space is zero-dimensional. For example, every topological space metrisable by an ultrametric is a first countable strongly zero-dimensional Hausdorff space.

\begin{crl}
Suppose that $X$ is strongly zero-dimensional. For $\mathscr{U}, \mathscr{U}' \in \m{UF}(|X|)$, the equality $\lim_{\mathscr{U}} |f(x)| = \lim_{\mathscr{U}'} |f(x)|$ holds for any $f \in \m{C}_{\m{bd}}(X,k)$ if and only if $F \cap F' \neq \emptyset$ for any closed subsets $F, F ' \subset X$ with $F \in \mathscr{U}$ and $F' \in \mathscr{U}'$.
\end{crl}

This is a generalisation of \cite{EM1} Theorem 4 for the class of bounded continuous functions, and together with Corollary \ref{uf - suf} implies \cite{EM1} Theorem 1. We remark if we removed the assumption of the strong zero-dimensionality, then there are obvious counter-examples. For example, a connected space is never strongly zero-dimensional unless it has at most one point, and every $k$-valued continuous function on a connected space is a constant function. In particular, every set-theoretical ultrafilter gives the same limit.

\begin{proof}
To begin with, suppose that the equality $\lim_{\mathscr{U}} |f(x)| = \lim_{\mathscr{U}'} |f(x)|$ holds for any $f \in \m{C}_{\m{bd}}(X,k)$. Then we have $\mathscr{U} \cap \m{CO}(X) = \mathscr{U}' \cap \m{CO}(X)$ by Corollary \ref{uf - suf}. Let $F, F ' \subset X$ be closed subsets with $F \in \mathscr{U}$ and $F' \in \mathscr{U}'$. Assume $F \cap F' = \emptyset$. Then there is a $U \in \m{CO}(X)$ such that $F \subset U \subset X \backslash F'$ because $X$ is strongly zero-dimensional. We obtain $U \in \mathscr{U} \cap \m{CO}(X)$ and $X \backslash U \in \mathscr{U}' \cap \m{CO}(X)$, and hence 
\begin{eqnarray*}
  \lim_{\mathscr{U}} |1_U(x)| = \| 1_U \|_{\mathscr{U} \cap \m{CO}(X)} = 1 \neq 0 = \| 1_U \|_{\mathscr{U}' \cap \m{CO}(X)} = \lim_{\mathscr{U}'} |1_U(x)|,
\end{eqnarray*}
where $1_U \colon X \to k$ denotes the characteristic function of $U$. It contradicts the assumption. Thus $F \cap F' \neq \emptyset$.

\vspace{0.2in}
Next, suppose that $F \cap F' \neq \emptyset$ for any closed subsets $F, F ' \subset X$ with $F \in \mathscr{U}$ and $F' \in \mathscr{U}'$. In order to verify $\mathscr{U} \cap \m{CO}(X) = \mathscr{U}' \cap \m{CO}(X)$, it suffices to show $\mathscr{U} \cap \m{CO}(X) \subset \mathscr{U}' \cap \m{CO}(X)$ by symmetry. Let $U \in \mathscr{U} \cap \m{CO}(X)$. Since $U \cap (X \backslash U) = \emptyset$, we have $X \backslash U \notin \mathscr{U}' \cap \m{CO}(X)$ by the assumption. Therefore $U \in \mathscr{U}' \cap \m{CO}(X)$ by the maximality of an ultrafilter. Thus $\mathscr{U} \cap \m{CO}(X) \subset \mathscr{U}' \cap \m{CO}(X)$, and hence $\mathscr{U} \cap \m{CO}(X) = \mathscr{U}' \cap \m{CO}(X)$. It implies that the equality $\lim_{\mathscr{U}} |f(x)| = \lim_{\mathscr{U}'} |f(x)|$ holds for any $f \in \m{C}_{\m{bd}}(X,k)$ by Corollary \ref{uf - suf}.
\end{proof}

\begin{crl}
The residue field of a maximal ideal of $\m{C}_{\m{bd}}(X,k)$ is $k$ if and only if it is a finite extension of $k$.
\end{crl}

This is a generalisation of \cite{EM2} Theorem 3.7 for the class of bounded continuous functions.

\begin{proof}
Let $m \subset \m{C}_{\m{bd}}(X,k)$ be a maximal ideal whose residue field is a finite extension $K$ of $k$. Take an arbitrary $f \in K$. Since $K$ is a finite extension of $k$, $f$ is algebraic over $k$. We prove $f \in k$. Assume $f \notin k$. Let $P(T) \in k[T]$ denote the minimal polynomial of $f$ over $k$. Let $L$ denote a decomposition field $P$, and fix an embedding $K \hookrightarrow L$. We endow $L$ with a unique extension of the valuation of $k$. Since $f \notin k$, $P(T)$ is an irreducible polynomial over $k$ with zeros $f_1, \ldots, f_d$ in $L \backslash k$. Since $L$ is a finite extension of $k$, $k$ is closed in $L$. Therefore for any $i \in \N \cap [1,d]$, the map $\xi_i \colon k \mapsto [0,\infty) \colon a \mapsto | a - f_i |$ is a continuous map with $\xi_i(a) \geq r_i$ for any $a \in k$ for some $r_i \in (0,\infty)$. In particular, we have $|P(a)| = \prod_{i = 1}^{d} \xi_i(a) \geq \prod_{i = 1}^{d} r_i > 0$. On the other hand, since $K$ is the residue field of $m$, there is an $F \in \m{C}_{\m{bd}}(X,k)$ whose image in $K$ is $f$. Then $F$ satisfies $P(F) \in m$. By the proof of Proposition \ref{maximal and utrafilter}, $m$ coincides with the support of the bounded multiplicative seminorm $\| \cdot \|_{\m{Ch}_m}$, and hence there is a $U \in \m{Ch}_m$ such that $\sup_{x \in U} |P(F)(x)| < \prod_{i = 1}^{d} r_i$ by the definition of $\| \cdot \|_{\m{Ch}_m}$. Since $U \neq \emptyset$, there exists an $x \in U$. However, we have $F(x) \in k$, and hence $|P(F)(x)| = |P(F(x))| \geq \prod_{i = 1}^{d} r_i$. It is a contradiction. Thus $f \in k$. We conclude that $K = k$.
\end{proof}

\begin{crl}
\label{partition of unity}
An ideal $I \subset \m{C}_{\m{bd}}(X,k)$ coincides with $\m{C}_{\m{bd}}(X,k)$ if and only if $I$ satisfies
\begin{eqnarray*}
  \inf_{x \in X} \sup_{f \in S} |f(x)| > 0
\end{eqnarray*}
for some non-empty finite subset $S \subset I$.
\end{crl}

This is a generalisation of \cite{EM1} Theorem 5 for the class of bounded continuous functions.

\begin{proof}
The sufficient implication is obvious because $1 \in \m{C}_{\m{bd}}(X,k)$. Suppose that $I$ does not coincides with $\m{C}_{\m{bd}}(X,k)$. Take a maximal ideal $m \subset \m{C}_{\m{bd}}(X,k)$ containing $I$. Let $S \subset m$ be a finite subset. Since $\| \cdot \|_{\m{Ch}_m}$ satisfies $\| f \|_{\m{Ch}_m} = 0$ for any $f \in m$, we have that for any $\epsilon \in (0,\infty)$, there is a $U \in \m{Ch}_m$ such that $\sup_{x \in U} |f(x)| < \epsilon$ for any $f \in S$. In particular, we obtain $\inf_{x \in X} \sup_{f \in S} |f(x)| = 0$ for any non-empty finite subset $S \subset I$.
\end{proof}

We remark that Corollary \ref{partition of unity} is also verified in a direct way with no use of our results. Indeed, let $I \subset \m{C}_{\m{bd}}(X,k)$ be an ideal such that there is a non-empty finite subset $S \subset I$ with $r \coloneqq \inf_{x \in X} \sup_{f \in S} |f(x)| > 0$. We put $U_f \coloneqq \Set{x \in X}{|f(x)| \geq r} \in \m{CO}(X)$ for each $f \in S$. Then by the assumption, the family $\mathscr{U} \coloneqq \set{U_f}{f \in S}$ covers $X$. Taking a total order on $S$, we put $S = \{ f_0, \ldots, f_d \}$. Then setting $U_i \coloneqq U_{f_i} \backslash \bigcup_{j = 0}^{i-1} U_{f_j}$ for each $i \in \N \cap [0,d]$, we obtain a refinement $\{ U_0, \ldots, U_d \}$ of $\mathscr{U}$ consisting of pairwise disjoint clopen subsets. For each $i \in \N \cap [0,d]$, we have $|f(x)| \geq r$ for any $x \in U_i$, and hence $g_i \coloneqq (1-1_{U_i}) + 1_{U_i} f$ is an invertible element of $\m{C}_{\m{bd}}(X,k)$ with $\| g_i^{-1} \| \leq \max \{ r^{-1}, 1 \}$, where $1_{U_i} \colon X \to k$ denotes the characteristic function of $U_i$. We obtain
\begin{eqnarray*}
  1 = \sum_{i = 0}^{d} 1_{U_i} = \sum_{i = 0}^{d} 1_{U_i} g_i^{-1} f_i \in I,
\end{eqnarray*}
and thus $I = \m{C}_{\m{bd}}(X,k)$.

\section{Related Results}
\label{Related Results}

\subsection{Another Construction}
\label{Another Construction}

In the case where $k$ is a local field or a finite field, we show that $\m{BSC}_k(X)$ coincides with a space $\m{SC}_k(X)$ defined in this section. Here a {\it local field} means a complete valuation field with non-trivial discrete valuation and finite residue field.

\begin{dfn}
\label{SC}
Denote by $\m{C}_{\m{bd}}(X,k)(1) \subset \m{C}_{\m{bd}}(X,k)$ the subset $\m{C}(X,k^{\circ})$ of bounded continuous $k$-valued functions on $X$ which take values in the subring $k^{\circ} \subset k$ of integral elements, and consider the evaluation map
\begin{eqnarray*}
  \iota'_k \colon X &\to& (k^{\circ})^{\m{C}_{\m{bd}}(X,k)(1)}\\
  x &\mapsto& (f(x))_{f \in \m{C}_{\m{bd}}(X,k)(1)}.
\end{eqnarray*}
By the definition of the direct product topology, $\iota'_k$ is continuous. Denote by $\m{SC}_k(X) \subset (k^{\circ})^{\m{C}_{\m{bd}}(X,k)(1)}$ the closure of the image of $\iota'_k$. We also denote by $\iota'_k$ the continuous map $X \to \m{SC}_k(X)$ induced by $\iota'_k$.
\end{dfn}

If $k$ is a local field or a finite field, then $\m{SC}_k(X)$ is a totally disconnected compact Hausdorff space because so is $k^{\circ}$.

\begin{prp}
\label{extension property}
The space $\m{SC}_k(X)$ satisfies the following extension property: For any $f \in \m{C}_{\m{bd}}(X,k)$, there is a unique $\m{SC}_k(f) \in \m{C}_{\m{bd}}(\m{SC}_k(X),k)$ such that $f = \m{SC}_k(f) \circ \iota'_k$. Moreover, the equality $\| f \| = \| \m{SC}_k(f) \|$ holds.
\end{prp}

\begin{proof}
The uniqueness of $\m{SC}_k(f)$ and the norm-preserving property is obvious because $\iota'_k(X) \subset \m{SC}_k(X)$ is dense and $k$ is Hausdorff. We construct the extension $\m{SC}_k(f)$. Note that $|k| \subset \lbrack 0, \infty )$ is bounded if and only if $|k| = \{ 0,1 \}$. Therefore $|k| \subset \lbrack 0, \infty )$ is unbounded or closed. It implies that that there is an $a \in k^{\times}$ such that $\| f \| \leq |a|$. For an $x = (x_g)_{g \in \m{C}_{\m{bd}}(X,k)(1)} \in \m{SC}_k(X)$, the value $a x_{a^{-1}f} \in k$ is independent of the choice of an $a \in k^{\times}$, and we set $\m{SC}_k(f)(x) \coloneqq a x_{a^{-1}f}$. Indeed, let $a_1,a_2 \in k^{\times}$ and suppose $\| f \| \leq \min \{ |a_1|, |a_2| \}$. For any $y \in X$, one has $\iota'_k(y)_{a_1^{-1}f} = a_1^{-1}f(y)$ and $\iota'_k(y)_{a_2^{-1}f} = a_2^{-1}f(y)$. It implies $a_1 \iota'_k(y)_{a_1^{-1}f} = a_2 \iota'_k(y)_{a_2^{-1}f} \in k$. Since the image $\iota'_k(X) \subset \m{SC}_k(X)$ is dense, one obtains $a_1 x_{a_1^{-1}f} = a_2 x_{a_2^{-1}f} \in k$. By the discussion above, one gets $\m{SC}_k(f) \circ \iota'_k = f$. The map $\m{SC}_k(f)$ is continuous by the definition of $\m{SC}_k(X)$.
\end{proof}

\begin{crl}
\label{functoriality}
For a $(Y,\varphi) \in X/\m{Top}$, there is a unique continuous map
\begin{eqnarray*}
\begin{array}{cccc}
  \m{SC}_k(\varphi) \colon \m{SC}_k(X) \to \m{SC}_k(Y)
\end{array}
\end{eqnarray*}
such that $\m{SC}_k(\varphi) \circ \iota'_k = \iota'_k \circ \varphi$.
\end{crl}

\begin{proof}
The uniqueness of $\m{SC}_k(\varphi)$ follows from the facts that $X$ is dense in $\m{SC}_k(\varphi)$ and that $\m{SC}_k(Y)$ is Hausdorff. By Proposition \ref{extension property}, one has a unique continuous map
\begin{eqnarray*}
  \m{SC}_k(\varphi) \colon \m{SC}_k(X) \to (k^{\circ})^{\m{C}_{\m{bd}}(Y,k)(1)}
\end{eqnarray*}
extending the composite
\begin{eqnarray*}
\begin{array}{cccc}
  X \stackrel{\varphi}{\longrightarrow} Y \stackrel{\iota'_k}{\longrightarrow} \m{SC}_k(Y) \hookrightarrow (k^{\circ})^{\m{C}_{\m{bd}}(Y,k)(1)}.
\end{array}
\end{eqnarray*}
Its image lies in the closed subspace $\m{SC}_k(Y)$ because $X$ is dense in $\m{SC}_k(X)$. One obtains a continuous map $\m{SC}_k(\varphi) \colon \m{SC}_k(X) \to \m{SC}_k(Y)$ such that $\m{SC}_k(\varphi) \circ \iota'_k = \iota'_k \circ \varphi$.
\end{proof}

Thus one obtains a functor
\begin{eqnarray*}
\begin{array}{cccc}
  \m{SC}_k \colon & \m{Top} & \to & \m{TDCHTop} \\
  & Y & \rightsquigarrow & \m{SC}_k(Y) \\
  & (\varphi \colon Y \to Z) & \rightsquigarrow & \left( \m{SC}_k(\varphi) \colon \m{SC}_k(Y) \to \m{SC}_k(Z) \right)
\end{array}
\end{eqnarray*}
with an obvious natural transform $\iota_k \colon \m{id}_{\m{Top}} \to \m{SC}_k$. We compare $\m{BSC}_k$ with $\m{SC}_k$ in the case where $k$ is a local field or a finite field.

\begin{lmm}
\label{reflexivity 2}
Suppose that $k$ is a local field or a finite field endowed with the trivial valuation, and that $X$ is a totally disconnected compact Hausdorff space. Then $\iota'_k \colon X \to \m{SC}_k(X)$ is a homeomorphism.
\end{lmm}

\begin{proof}
By the assumption of $k$, $\m{SC}_k(X)$ is a totally disconnected compact Hausdorff space. Therefore it suffices to verify the injectivity of $\iota'_k$, because a continuous map from a compact space to a Hausdorff space is a closed map. Let $x,y \in X$ with $x \neq y$. Since $X$ is zero-dimensional and Hausdorff, there is a $U \in \m{CO}(X)$ such that $x \in U$ and $y \notin U$. Then one has $\iota'_k(x)_{1_U} = 1 \neq 0 = \iota'_k(y)_{1_U}$, and hence $\iota'_k(x) \neq \iota'_k(y)$. Thus $\iota'_k \colon X \to \m{SC}_k(X)$ is injective.
\end{proof}

\begin{prp}
\label{SC and TDC}
Suppose that $k$ is a local field or a finite field endowed with the trivial valuation. Then $\m{SC}_k(X)$ is initial in $X/\m{TDCHTop}$ with respect to $\iota'_k$.
\end{prp}

We remark that the assumption on the base field $k$ is not necessary when $X$ is compact. Analysis of continuous functions on a compact space is quite classical.

\begin{proof}
For a $(Y,\varphi) \in \m{ob}(X/\m{TDCHTop})$, we construct a continuous extension
\begin{eqnarray*}
  \psi \colon \m{SC}_k(X) \to Y
\end{eqnarray*}
of $\varphi$ in an explicit way. An extension $\psi$ is unique if it exists, because the image of $X$ is dense in $\m{SC}_k(X)$ and $Y$ is Hausdorff. Consider the commutative diagram
\begin{eqnarray*}
\begin{CD}
  X           @>{\varphi}>>           Y \\
  @V{\iota'_k}VV                      @VV{\iota'_k}V \\
  \m{SC}_k(X) @>>{\m{SC}_k(\varphi)}> \m{SC}_k(Y)
\end{CD}.
\end{eqnarray*}
By Lemma \ref{reflexivity 2}, the right vertical map is a homeomorphism, and one obtains a continuous map $\psi \coloneqq \iota'_k{}^{-1} \circ \m{SC}_k(\varphi) \colon \m{SC}_k(X) \to Y$.
\end{proof}

\begin{crl}
\label{SC and BSC}
Suppose that $k$ is a local field or a finite field endowed with the trivial valuation.
\begin{itemize}
\item[(i)] The space $\m{SC}_k(X)$ is homeomorphic to $\m{BSC}_k(X)$ under $X$.
\item[(ii)] The space $\m{BSC}_k(X)$ satisfies the extension property for a bounded continuous $k$-valued function on $X$ in Proposition \ref{extension property}.
\item[(iii)] The natural homomorphism $\m{C}(\m{BSC}_k(X),k) \to \m{C}_{\m{bd}}(X,k)$ is an isometric isomorphism.
\item[(iv)] The space $\m{BSC}_k(X)$ consists of $k$-rational points, and the residue field of any maximal ideal of $\m{C}_{\m{bd}}(X,k)$ is $k$.
\end{itemize}
\end{crl}

\begin{proof}
We deal only with (iv). Since every maximal ideal of $\m{C}_{\m{bd}}(X,k)$ is the support of an $x \in \m{BSC}_k(X)$ as is referred in the proof of Proposition \ref{image of supp}, it suffices to verify the first assertion. We recall that for a Banach $k$-algebra $\A$, an $x \in \M_k(\A)$ is said to be a $k$-rational point if its support $\set{f \in A}{x(f) = 0}$ is a maximal ideal of $\A$ whose residue field is $k$. The isomorphism $\m{C}(\m{BSC}_k(X),k) \to \m{C}_{\m{bd}}(X,k)$ in (iii) gives an identification $\m{BSC}_k(X) = \M_k(\m{C}(\m{BSC}_k(X),k))$. The assertion immediately follows from a non-Archimedean generalisation of Stone--Weierstrass theorem (\cite{Ber1} 9.2.5.\ Theorem (i)) for $\m{C}(\m{BSC}_k(X),k)$.
\end{proof}

In particular, concerning Corollary \ref{SC and BSC} (iv), any point of $\m{BSC}_k(X)$ is peaked in the sense that the complete residue field is a peaked Banach $k$-algebra (\cite{Ber1} 5.2.1,\ Definition). A Banach $k$-algebra $\A$ is said to be {\it peaked} if for any extension $K/k$ of complete valuation fields, the norm on $\A \hat{\otimes}_k K$ is multiplicative. The notion of a peaked point is useful when we consider the topology-theoretical multiplication of points of analytic group. Corollary \ref{SC and BSC} (iv) does not hold when $k = \C_p$. Indeed, consider the rigid analytic disc $\m{D}^1(\C_p) \coloneqq \M_{\C_p}(\C_p \{ z \})$. It admits a natural embedding $\C_p^{\circ} \hookrightarrow \m{D}^1(\C_p)$ into a dense subset. The bounded $\C_p$-algebra homomorphism $\varphi \colon \C_p \{ z \} \to \m{C}_{\m{bd}}(\C_p^{\circ},\C_p)$ sending the variable $z$ to the coordinate function $z \colon \C_p^{\circ} \hookrightarrow \C_p$ induces a continuous map $\varphi^* \colon \m{BSC}_{\C_p}(\C_p^{\circ}) \to \m{D}^1(\C_p)$ under $\C_p^{\circ}$. Since $\m{BSC}_{\C_p}(\C_p^{\circ})$ and $\m{D}^1(\C_p)$ are compact Hausdorff spaces, the image of $\varphi^*$ contains the closure of the dense subset $\C_p^{\circ} \subset \m{D}^1(\C_p)$. Therefore $\varphi^*$ is surjective. Since $\C_p$ is not spherically complete, there is a $y \in \m{D}^1(\C_p)$ of type $4$. Take an $x \in \m{BSC}_{\C_p}(\C_p^{\circ})$ in $\varphi^{-1}(y)$. The induced bounded $\C_p$-algebra homomorphism $\varphi_x \colon \C_p \{ z \}/\m{supp}(y) \to \m{C}_{\m{bd}}(\C_p^{\circ},\C_p)/\m{supp}(x)$ gives an extension of fields transcendental over $\C_p$ because $y$ is of type $4$. Thus the point $x \in \m{BSC}_{\C_p}(\C_p^{\circ})$ is not $\C_p$-rational.

\subsection{Relation to the Stone--$\check{\m{C}}$ech compactification}
\label{Relation to the Stone--Cech compactification}

A compact Hausdorff space is totally disconnected if and only if it is zero-dimensional, and hence under $X$, the notion of an initial totally disconnected compact Hausdorff space is equivalent to that of an initial zero-dimensional compact Hausdorff space. Banaschewski constructed a zero-dimensional compact Hausdorff space $\zeta X$ under $X$ in the case where $X$ is zero-dimensional in \cite{Ban}, and proved that $\zeta X$ is initial. In the case where $X$ is a zero-dimensional Hausdorff space, the structure map $X \to \zeta X$ is a homeomorphism onto the image and $\zeta X$ is sometimes called the Banaschewski compactification.

\vspace{0.1in}
Let $p$ be a prime number. By the argument above, $\m{BSC}_{\Q_p}(X)$ is one of the generalisation of $\zeta X$. The Stone--$\check{\m{C}}$ech compactification $\beta X$ has a universality as an initial object in the category of compact topological spaces under $X$, and hence there is a unique continuous map $\beta X \to \m{BSC}_{\Q_p}(X)$ under $X$. Then a natural question arises: ``When is the map $\beta X \to \m{BSC}_{\Q_p}(X)$ is a homeomorphism?'' In other words, ``When is $\beta X$ totally disconnected?'' In such a case, $\m{BSC}_{\Q_p}(X)$ satisfies the extension property for an Archimedean bounded continuous function on $X$. This connection between the Archimedean analysis and the non-Archimedean analysis looks interesting. It is well-known that if $X$ is an infinite discrete space, then $\beta X$ is totally disconnected. In particular, one has $\beta \N \cong \m{BSC}_{\Q_p}(\N)$. Furthermore, Banaschewski proved $\beta X$ is homeomorphic to $\zeta X$ if $X$ is a second countable zero-dimensional Hausdorff space (\cite{Ban} Satz 6.\ Korollar 2). In this case, $\beta X$ is totally disconnected and hence is homeomorphic to $\m{BSC}_{\Q_p}(X)$. For example, the closed unit disc $\C_{\ell}^{\circ} \subset \C_{\ell}$ for a prime number $\ell \in \N$ is a second countable zero-dimensional Hausdorff space, and hence one has $\beta \C_{\ell}^{\circ} \cong \m{BSC}_{\Q_p}(\C_{\ell}^{\circ})$.

\vspace{0.1in}
We do not know whether there is an example of a totally disconnected or zero-dimensional space $X$ such that the map $\beta X \to \m{BSC}_{\Q_p}(X)$ is not a homeomorphism. Banaschewski gave the following necessary condition for the bijectivity of $\beta X \to \m{BSC}_{\Q_p}(X)$ in \cite{Ban} Satz 2. For a normal zero-dimensional space $X$, if the continuous map $\beta X \to \m{BSC}_{\Q_p}(X)$ is a homeomorphism, then $X$ is of $\check{\m{C}}$ech-dimension $0$, i.e.\ any finite open covering of $X$ admits a finite clopen refinement. Therefore if there is a normal zero-dimensional space of positive $\check{\m{C}}$ech-dimension, then it is an example.

\section{Applications}
\label{Applications}

\subsection{Automatic Continuity Theorem}
\label{Automatic Continuity Theorem}

One of important classical problems for a commutative $C^*$-algebra is Kaplansky Conjecture on automatic continuity problem, which claims that every injective $\C$-algebra homomorphism $\varphi \colon \m{C}_{\m{bd}}(X,\C) \hookrightarrow \A$ is continuous for any Banach $\C$-algebra $\A$. Consider the following weak version: every injective $\C$-algebra homomorphism $\varphi \colon \m{C}_{\m{bd}}(X,\C) \hookrightarrow \A$ with closed image is continuous for any Banach $\C$-algebra $\A$. It was proved by Kaplansky, and Solovay showed that the conjecture is independent of the axiom of ZFC. For more details about the automatic continuity problem, see \cite{Dal}. \cite{Sol}, and \cite{Woo}. Now it is natural to consider an analogous question in the non-Archimedean case, and we prove its weak version in this subsection.

\begin{dfn}
Let $\A$ be a Banach $k$-algebra, and $m \subset \A$ a maximal ideal whose residue field $K$ is a finite extension of $k$. For each $f \in \A$, we denote by $f(m)$ the image of $f$ in $K$, and by $|f(m)|$ the norm of $f(m)$ with respect to a unique extension of the norm of $k$.
\end{dfn}

\begin{lmm}
\label{orthogonality of a maximal ideal of codimension 1}
Let $\A$ be a Banach $k$-algebra. For a maximal ideal $m \subset \A$ whose residue field is $k$, the canonical projection $\A \twoheadrightarrow \A/m \cong_k k$ gives a decomposition $\A = k \oplus m$ as the orthogonal direct sum, i.e.\ the equality
\begin{eqnarray*}
  \| a + g \| = \max \{ |a|, \| g \| \}
\end{eqnarray*}
holds for any $a \in k$ and $g \in m$.
\end{lmm}

\begin{proof}
Since the composite $k \hookrightarrow \A \twoheadrightarrow \A/m$ is a bijective $k$-linear homomorphism, one obtains a decomposition $\A = k \oplus m$ as the direct sum of $k$-vector spaces. Take an element $f \in \A$, and denote by $f(m) \in k$ the image of $f$ in the quotient $\A/m \cong_k k$. In order to prove the orthogonality of the direct sum, it suffices to show $\| f \| = \max \{ |f(m)|, \| f - f(m) \| \}$. The inequality $\leq$ is obvious. If $|f(m)| \neq \| f - f(m) \|$, the equality follows from the general property of a non-Archimedean norm, and hence we may assume $|f(m)| = \| f - f(m) \|$. If $f(m) = 0$, then one has $\| f - f(m) \| = 0$ and therefore $f = f(m) + (f - f(m)) = 0$. Suppose $f(m) \neq 0$. Assume $\| f \| < |f(m)|$. Then one has $\| f(m)^{-1}f \| < 1$, and hence
\begin{eqnarray*}
  f - f(m) = -f(m)(1 - f(m)^{-1}f) \in k^{\times} \A^{\times} = \A^{\times}
\end{eqnarray*}
by \cite{BGR} 1.2.4.\ Proposition 4. It contradicts the fact $f - f(m) \in m$, and thus $\| f \| = |f(m)| = \max \{ |f(m)|, \| f - f(m) \| \}$.
\end{proof}

We may apply Lemma \ref{orthogonality of a maximal ideal of codimension 1} to any maximal ideal of $\m{C}_{\m{bd}}(X,k)$ by Corollary \ref{SC and BSC} (iv).

\begin{crl}
\label{orthogonal}
Suppose that $k$ is a local field or a finite field endowed with the trivial norm. For any maximal ideal $m \subset \m{C}_{\m{bd}}(X,k)$, the canonical projection $\m{C}_{\m{bd}}(X,k) \twoheadrightarrow \m{C}_{\m{bd}}(X,k)/m$ gives a decomposition $\m{C}_{\m{bd}}(X,k) = k \oplus m$ as the orthogonal direct sum of $k$-Banach spaces.
\end{crl}

This is a partial generalisation of \cite{EM1} Theorem 7, which states that if $X$ is an ultrametric space and $k$ is locally compact, then the same holds.

\begin{prp}
\label{algebraic BFA}
Suppose that $k$ is a local field or a finite field endowed with the trivial norm. Let $\m{Max}(\m{C}_{\m{bd}}(X,k)) \subset \m{Spec}(\m{C}_{\m{bd}}(X,k))$ denote the subset of maximal ideals. Then the equality
\begin{eqnarray*}
  \| f \| = \sup_{m \in \m{Max}(\m{C}_{\m{bd}}(X,k))} |f(m)|
\end{eqnarray*}
holds for any $f \in \m{C}_{\m{bd}}(X,k)$. In particular, the norm of $\m{C}_{\m{bd}}(X,k)$ is determined by the algebraic structure of it.
\end{prp}

\begin{proof}
Since the norm of $\m{C}_{\m{bd}}(X,k)$ is power-multiplicative, the equality
\begin{eqnarray*}
  \| f \| = \sup_{x \in \m{BSC}(X)} x(f)
\end{eqnarray*}
holds for any $f \in \m{C}_{\m{bd}}(X,k)$ by \cite{Ber1} 1.3.1.\ Theorem. This gives the assertion because $\m{supp} \colon \m{BSC}_k(X) \to \m{Spec}(\m{C}_{\m{bd}}(X,k))$ is bijective onto $\m{Max}(\m{C}_{\m{bd}}(X,k))$ by the proof of Theorem \ref{main theorem}, and because the condition (iv) in the definition of a bounded multiplicative seminorm in \S \ref{Berkovich Spectra} guarantees $x(f) = |f(\m{supp}(x))|$ by Corollary \ref{SC and BSC} (iv).
\end{proof}

\begin{prp}
\label{equivalence of complete norms}
Suppose that $k$ is a local field. Then every complete norm on the underlying $k$-algebra of $\m{C}_{\m{bd}}(X,k)$ is equivalent to each other.
\end{prp}

\begin{proof}
Since $k$ is a local field, the norm of $k$ is not trivial, and hence the boundedness of a $k$-linear homomorphism between normed $k$-vector spaces is equivalent to the continuity. Therefore it suffices to show that the identity $\m{id} \colon \m{C}_{\m{bd}}(X,k) \to \m{C}_{\m{bd}}(X,k)$ is a homeomorphism with respect to the metric topologies given by an arbitrary complete norm $\| \cdot \|'$ of the source and the supremum norm $\| \cdot \|$ of the target. By Lemma \ref {orthogonality of a maximal ideal of codimension 1} applied to $(\m{C}_{\m{bd}}(X,k),\| \cdot \|')$ and Corollary \ref{algebraic BFA}, $\m{id}$ is a $k$-linear contraction map, and hence is continuous. Moreover, since the norm of $k$ is not trivial, the open mapping theorem holds by \cite{BGR} 2.8.1.\ Theorem, and therefore $\m{id}$ is an open map. Thus $\m{id}$ is a homeomorphism.
\end{proof}

\begin{thm}
\label{automatic continuity}
Suppose that $k$ is a local field, and let $\A$ be a Banach $k$-algebra. Then any injective $k$-algebra homomorphism $\varphi \colon \m{C}_{\m{bd}}(X,k) \hookrightarrow \A$ whose image is closed is continuous.
\end{thm}

\begin{proof}
Since the underlying metric spaces of $\m{C}_{\m{bd}}(X,k)$ and $\A$ are complete, it suffices to show that $\varphi$ sends a Cauchy sequence in $\m{C}_{\m{bd}}(X,k)$ to a Cauchy sequence in $\A$. Let $\| \cdot \|' \colon \m{C}_{\m{bd}}(X,k) \to [0,\infty)$ denote the composite of $\varphi$ and the norm of $\A$. Then since $\varphi$ is a bijective homomorphism of $k$-algebras, $\| \cdot \|'$ is a norm of the $k$-algebra $\m{C}_{\m{bd}}(X,k)$. Since the image of $\varphi$ is closed, $\| \cdot \|'$ is complete, and hence is equivalent to the supremum norm by Proposition \ref{equivalence of complete norms}. We conclude that $\varphi$ is continuous.
\end{proof}

The automatic continuity theorem immediately yields a criterion for the continuity of a faithful representation over a local field.

\begin{crl}
\label{automatic continuity of Banach representations}
Suppose that $k$ is a local field. Let $V$ be a Banach $k$-vector space and $\rho \colon \m{C}_{\m{bd}}(X,k) \times V \to V$ a $k$-linear representation of a $k$-algebra $\m{C}_{\m{bd}}(X,k)$ satisfying the following conditions:
\begin{itemize}
\item[(i)] The $k$-linear operator $\tilde{\rho}_f \colon V \to V \colon v \mapsto \rho(f,v)$ is bounded for any $f \in \m{C}_{\m{bd}}(X,k)$.
\item[(ii)] The $k$-linear representation $\rho$ is faithful, i.e.\ the equality $\rho_f = 0$ implies $f = 0$ for any $f \in \m{C}_{\m{bd}}(X,k)$.
\item[(iii)] The image of the induced $k$-algebra homomorphism $\tilde{\rho}_{\cdot} \colon \m{C}_{\m{bd}}(X,k) \to \B_k(V)$ is closed, where $\B_k(V)$ is the Banach $k$-algebra of bounded operators on $V$.
\end{itemize}
Then $\tilde{\rho}_{\cdot}$ is bounded, and in particular, $\rho \colon \m{C}_{\m{bd}}(X,k) \times V \to V$ is continuous.
\end{crl}

\subsection{Ground Field Extensions}
\label{Ground Field Extensions}

We study the ground field extensions of $\m{C}_{\m{bd}}(X,k)$. We note that there are two distinct notions of the ground field extensions. One is given by extending the scalar of functions, and the other is given by tensoring the scalar.

\begin{prp}
\label{independent of k}
Let $K$ and $L$ be complete valuation fields. Then there exists a unique homeomorphism $\m{BSC}_K(X) \cong \m{BSC}_L(X)$ compatible with the evaluation maps.
\end{prp}

We remark that we do not assume that $K$ and $L$ contains the same base field $k$, and hence, for example, it is possible to choose $\Q_p$ and $\F_{\ell}((T))$ for $K$ and $L$ respectively.

\begin{proof}
The assertion holds because $\m{BSC}_K(X)$ and $\m{BSC}_L(X)$ are initial objects with repsect to the evaluation maps in $X/\m{TDCHTop}$ by Corollary \ref{initial 2}.
\end{proof}

\begin{prp}
\label{the extension of the scalar}
Let $K/k$ be an extension of complete valuation fields. Then the ground field extension $\m{BSC}_K(X) \to \m{BSC}_k(X)$ associated with the natural embedding $\m{C}_{\m{bd}}(X,k) \hookrightarrow \m{C}_{\m{bd}}(X,K)$ is a homeomorphism.
\end{prp}

\begin{proof}
The ground field extension above is compatible with the evaluation maps, and coincides with the unique homeomorphism in Proposition \ref{independent of k}.
\end{proof}

Now we consider the other ground field extension, namely, the canonical $K$-algebra homomorphism $K \hat{\otimes}_k \m{C}_{\m{bd}}(X,k) \to \m{C}_{\m{bd}}(X,K)$ induced by the universal property of the complete tensor product in the category of Banach $k$-algebras. In fact, the ground field extension is not an isomorphism in general, and it yields a criterion for a topological property of $X$ and the valuation of $k$.

\begin{lmm}
\label{locally constant}
The $k$-subalgebra of $\m{C}_{\m{bd}}(X,k)$ consisting of locally constant bounded functions is dense.
\end{lmm}

\begin{proof}
Take an $f \in \m{C}_{\m{bd}}(X,k)$. If $f = 0$, then $f$ is locally constant. Suppose $f \neq 0$. For an $\epsilon > 0$, the pre-image of every open disc of radius $\epsilon$ in $k$ by $f$ is clopen. Therefore one obtains a pairwise disjoint clopen covering $\mathscr{U}$ of $X$ such that the image $f(U)$ is contained in an open disc of radius $\epsilon$ in $k$ for any $U \in \mathscr{U}$. Fix an $a_U \in f(U)$ for each $U \in \mathscr{U}$. The infinite sum $g \coloneqq \sum_{U \in \mathscr{U}} a_U1_U \colon X \to k$ convergences pointwise to a locally constant continuous function.  The obvious inequality $|g(x)| \leq \| f \|$ holds for any $x \in X$, and hence $g$ is bounded. One has $\| f - g \| \leq \epsilon$ by the definition of the disjoint clopen covering $\mathscr{U}$, and hence the $k$-subalgebra of locally constant functions is dense in $\m{C}_{\m{bd}}(X,k)$.
\end{proof}

\begin{lmm}
\label{two ground field extensions}
Suppose that $k$ is spherically complete (\cite{BGR} 2.4.4.\ Definition 1). Let $K/k$ be an extension of complete valuation fields. Then the natural bounded $K$-algebra homomorphism $\iota_{K/k} \colon K \hat{\otimes}_k \m{C}_{\m{bd}}(X,k) \to \m{C}_{\m{bd}}(X,K)$ is an isometry.
\end{lmm}

For example, a local field and every field endowed with the trivial norm are spherically complete. We will use this lemma for $\F_p$ endowed with the trivial norm, $\Q$ endowed with the trivial norm, and $\Q_p$.

\begin{proof}
Take an $f \in K \hat{\otimes}_k \m{C}_{\m{bd}}(X,k)$. If $f = 0$, then $\| \iota_{K/k}(f) \| = 0 = \| f \|$, and hence we assume $f \neq 0$. In particular, $X \neq \emptyset$ and both of $K \hat{\otimes}_k \m{C}_{\m{bd}}(X,k)$ and $\m{C}_{\m{bd}}(X,K)$ are non-zero Banach $K$-algebras. Therefore the norm of the bounded $K$-algebra homomorphism $\iota_{K/k}$ is $1$ because $\iota_{K/k}(1) = 1$ and the power-multiplicativity of the norm of $\m{C}_{\m{bd}}(X,K)$ guaranties that $\iota_{K/k}$ is submetric. Set $\epsilon \coloneqq \| f \|/2$, and take an element $g = \sum_{i = 1}^{n} a_i \otimes g_i \in K \otimes_k \m{C}_{\m{bd}}(X,k)$ with $\| f - g \| < \epsilon$ in $K \hat{\otimes}_k \m{C}_{\m{bd}}(X,k)$. We may assume $a_i \neq 0$ for any $i = 1,\ldots,n$ without loss of generality. By Lemma \ref{locally constant}, there is a locally constant bounded $k$-valued function $g'_i \in \m{C}_{\m{bd}}(X,k)$ such that $\| g_i - g'_i \| < |a_i|^{-1} \epsilon$ for each $i = 1,\ldots,n$. In particular, setting $g' \coloneqq \sum_{i = 1}^{n} a_i \otimes g'_i \in K \otimes_k \m{C}_{\m{bd}}(X,k)$, one has
\begin{eqnarray*}
  && \| f - g' \| = \| (f - g) + (g - g') \| \leq \max \left\{ \| f - g \|, \left\| \sum_{i = 1}^{n} a_i(g_i - g'_i) \right\| \right\}\\
  &\leq& \max \{ \| f - g \|, \max{}_{i = 1}^{n} |a_i| \| g_i - g'_i \| \} < \epsilon < \| f \|,
\end{eqnarray*}
and hence $\| g' \| = \| f \|$. Since $k$ is spherically complete, the finite dimensional normed $k$-vector subspace $k a_1 + \cdots + k a_n \subset K$ is $k$-Cartesian by \cite{BGR} 2.4.4.\ Proposition 2, and hence there is an orthogonal basis $b_1,\ldots,b_m \in K$ of $k a_1 + \cdots + k a_n$. Expressing $a_1,\ldots,a_n$ as a $k$-linear combination of $b_1,\ldots,b_m$, one obtains an expression $g' = \sum_{i = 1}^{m} b_ig''_i$ by a unique system $g''_1,\ldots,g''_m \in \m{C}_{\m{bd}}(X,k)$ of $k$-valued locally constant functions. For any $x \in X$, one has
\begin{eqnarray*}
  |\iota_{K/k}(g')(x)| = \left| \sum_{i = 1}^{m} b_i g''_i(x) \right| = \max{}_{i = 1}^{m} |g''_i(x)| |b_i|
\end{eqnarray*}
by the orthogonality of $b_1,\ldots,b_m$, and hence
\begin{eqnarray*}
  && \| \iota_{K/k}(g') \| = \sup_{x \in X} |\iota_{K/k}(g')(x)| = \sup_{x \in X} \max{}_{i = 1}^{m} |b_i| |g''_i(x)| = \max{}_{i = 1}^{m} |b_i| \sup_{x \in X} |g''_i(x)|\\
  &=& \max{}_{i = 1}^{m} |b_i| \|g''_i \| \geq \| g' \|.
\end{eqnarray*}
Since $\iota_{K/k}$ is a $K$-linear contraction map, one gets $\| \iota_{K/k}(g') \| = \| g' \|$. We conclude
\begin{eqnarray*}
  \| \iota_{K/k}(f - g') \| \leq \| f - g' \| < \epsilon < \| f \| = \| g' \| = \| \iota_{K/k}(g') \|
\end{eqnarray*}
and thus
\begin{eqnarray*}
  \| \iota_{K/k}(f) \| = \| \iota_{K/k}(f - g') + \iota_{K/k}(g') \| = \| \iota_{K/k}(g') \| = \| g' \| = \| f \|.
\end{eqnarray*}
\end{proof}

We denote by $\F \subset k$ the topological closure of the field $F$ of fractions of the image of the canonical ring homomorphism $\Z \to k$. We remark that $F$ is $\F_p$ if and only if $k$ is of characteristic $p > 0$, and is $\Q$ otherwise. In the former case, $\F$ is $\F_p$ endowed with the trivial valuation. In the latter case, $\F$ is $\Q$ endowed with the trivial norm if and only if $k$ is of equal characteristic $(0,0)$, and is $\Q_p$ if and only if $k$ is of mixed characteristic $(0,p)$. In particular, $\F$ is spherically complete. We determine when $\iota_{k/\F} \colon k \hat{\otimes}_{\F} \m{C}_{\m{bd}}(X,\F) \to \m{C}_{\m{bd}}(X,k)$ is an isomorphism. The following shows that $\m{C}_{\m{bd}}(X,k)$ is ``naive'' enough if and only if $X$ is compact or $k$ is sufficiently small in some sense.

\begin{thm}
\label{main theorem 2}
Suppose that $X$ is zero-dimensional and Hausdorff. Then the following are equivalent:
\begin{itemize}
\item[(i)] The space $X$ is compact, or $k$ is a local field or a finite field endowed with the trivial norm.
\item[(ii)] The $k$-subalgebra of $\m{C}_{\m{bd}}(X,k)$ generated by idempotents is dense.
\end{itemize}
In addition if $\F \neq \Q$, then (i) and (ii) are equivalent to the following:
\begin{itemize}
\item[(iii)] The homomorphism $\iota_{k/\F} \colon k \hat{\otimes}_{\F} \m{C}_{\m{bd}}(X,\F) \to \m{C}_{\m{bd}}(X,k)$ is an isometric isomorphism.
\item[(iv)] The space $\m{BSC}_k(X)$ consists of $k$-rational points.
\item[(v)] The map $\iota_k \colon X \hookrightarrow \m{BSC}_k(X)$ induces an isometric isomorphism $\m{C}(\m{BSC}_k(X),k) \to \m{C}_{\m{bd}}(X,k)$.
\item[(vi)] The space $\m{BSC}_k(X)$ satisfies the extension property for a bounded continuous $k$-valued functions in Proposition \ref{extension property}.
\end{itemize}
\end{thm}

We remark that a similar relation between (i) and (iv) is also verified by Alain Escassut and Nicolas Ma\"inetti in \cite{EM1} in the case where $X$ is an ultrametric space. For example, they proved in Theorem 7 that if $X$ is an ultrametric space and $k$ is locally compact, then (iv) holds.

\begin{proof}
Assume (i). We verify (ii). Take an $f \in \m{C}_{\m{bd}}(X,k)$. If $k$ is a local field or a finite field, the closed disc $\set{a \in k}{|a| \leq \| f \|} \subset k$ is compact. Otherwise $X$ is compact. Therefore for an $\epsilon > 0$, there is a finite pairwise disjoint clopen covering $\mathscr{U}$ of $X$ such that the image $f(U) \subset k$ is contained in an open disc of radius $\epsilon$ in $k$ for any $U \in \mathscr{U}$. Fixing an $a_U \in f(U)$ for each $U \in \mathscr{U}$, one has $\| f - \sum_{U \in \mathscr{U}} a_U1_U \| < \epsilon$. Thus  $k$-subalgebra of $\m{C}_{\m{bd}}(X,k)$ generated by idempotents is dense.

\vspace{0.1in}
Assume (ii). We verify (iii). Since $\iota_{k/\F}$ is an isometry by Lemma \ref{two ground field extensions}, it suffices to show that the image of $\iota_{k/\F}$ is dense. An idempotent of $\m{C}_{\m{bd}}(X,k)$, which is a characteristic function on a clopen subset of $X$, is contained in the subset $\m{C}_{\m{bd}}(X,\F) \subset \m{C}_{\m{bd}}(X,k)$. Therefore the image of the natural homomorphism $k \otimes_{\F} \m{C}_{\m{bd}}(X,\F) \to \m{C}_{\m{bd}}(X,k)$ is dense by (ii), and hence the image of $\iota_{k/\F}$ is dense.

\vspace{0.1in}
Assume (iii). We verify (iv) in the case $\F \neq \Q$. For an $x \in \m{BSC}_k(X)$, consider the composite $x' \colon \m{C}_{\m{bd}}(X,\F) \to k(x)$ of $x$ and the natural embedding $\m{C}_{\m{bd}}(X,\F) \hookrightarrow \m{C}_{\m{bd}}(X,k)$, where $k(x)$ is the completed residue field at $x$. Since $\F$ is contained in $k$, the character $x'$ defines an element $x' \in \m{BSC}_{\F}(X)$. Recall that $\F = \F_p$ or $\Q_p$ now. Since $\m{BSC}_{\F}(X)$ consists of $\F$-rational points by Corollary \ref{SC and BSC} (iv), the image of $x'$ is contained in $\F \subset k$. Therefore (iii) guarantees that the image of $x$ is contained in the closure of the $k$-vector subspace of $k(x)$ generated by $\F \subset k$, namely, the $1$-dimensional vector subspace $k \subset k(x)$. It implies $k(x) = k$.

\vspace{0.1in}
Assume (iv). We verify (v) in the case $\F \neq \Q$, and hence suppose that $\m{BSC}_k(X)$ consists of $k$-rational points. Then the Gel'fand transform $\m{C}_{\m{bd}}(X,k) \to \m{C}(\m{BSC}_k(X),k)$ is an isometric isomorphism by \cite{Ber1} 9.2.7.\ Corollary (ii), and coincides with the bounded $k$-algebra homomorphism induced by $\iota_k$.

\vspace{0.1in}
Assume (v). We verify (vi) in the case $\F \neq \Q$, and hence suppose that $\iota_k$ induces an isometric isomorphism $\m{C}(\m{BSC}_k(X),k) \to \m{C}_{\m{bd}}(X,k)$. Take an $f \in \m{C}_{\m{bd}}(X,k)$. The extension of $f$ on $\m{BSC}_k(X)$ is unique because the image of $\iota_k$ is dense by Corollary \ref{dense} and $k$ is Hausdorff. Since $\iota_k$ induces an isomorphism $\m{C}(\m{BSC}_k(X),k) \to \m{C}_{\m{bd}}(X,k)$, there is an $f' \in \m{C}(\m{BSC}_k(X),k)$ whose image is $f$, or in other words, $f'$ is the extension of $f$ on $\m{BSC}_k(X)$.

\vspace{0.1in}
Assume that (vi) and $\F \neq \Q$ hold, or that (ii) and $\F = \Q$ hold. We verify (i). Suppose that $X$ is non-compact and $k$ is neither a local field nor a finite field. Since $X$ is a zero-dimensional and non-compact, there is an $\mathscr{F} \in \m{UF}(X)$ without a cluster point by Proposition \ref{ultrafilter criterion}. In particular, $\mathscr{F}$ contains an infinite descending chain $X = U_0 \supsetneq U_1 \supsetneq \cdots$. Indeed, for any $U \in \mathscr{F}$ and $x \in U \neq \emptyset$, there is a $V \in \m{CO}(X)$ such that $V \in \mathscr{F}$, $V \subset U$, and $x \notin V$ because $x$ is not a cluster point of $\mathscr{F}$. One obtains an infinite set $\mathscr{U} = \set{U_i \backslash U_{i + 1}}{i \in \N}$ of pairwise disjoint clopen subsets of $X$. If the residue field $\tilde{k}$ of $k$ is an infinite field, set $Y \coloneqq \tilde{k}$ and take a set-theoretical lift $\varphi \colon Y \hookrightarrow k^{\circ}$ of the canonical projection $k^{\circ} \twoheadrightarrow Y$. Otherwise, the image $|k^{\times}| \subset (0,\infty)$ is dense because $k$ is neither a local field nor a finite field. Set $Y \coloneqq |k^{\times}| \cap (1/2,1) \subset (0,\infty)$, and take a set-theoretical lift $\varphi \colon Y \hookrightarrow k^{\circ}$ of the norm $|\cdot| \colon k \to \lbrack 0,\infty )$. Since $Y$ is dense in $(1/2,1)$, it is an infinite set. In both cases, endow $Y$ with the discrete topology. Since $Y$ is an infinite set, there is an injective map $\psi \colon \N \hookrightarrow Y$. The composite $\varphi \circ \psi \colon \N \hookrightarrow k^{\circ}$ is an injective continuous map, and the image is a closed discrete subspace because $|\varphi(y) - \varphi(y')| > 1/2$ for any $y,y' \in Y$. Since $\mathscr{U} \subset \m{CO}(X)$ is an infinite covering of $X$, there is an injective map $\Psi \colon \N \hookrightarrow \mathscr{U}$. Then the pointwise convergent infinite sum
\begin{eqnarray*}
  f \coloneqq \sum_{n \in \N} \varphi(\psi(n))1_{\Psi(n)} \colon X \to k
\end{eqnarray*}
determines a locally constant bounded function on $X$. There is a non-principal ultrafilter $\mathscr{F} \in \m{CO}(\N)$ by Lemma \ref{non-principal ultrafilter}. Now assume $\F \neq \Q$. By the conditions (vi), there is a continuous extension $\m{BSC}_k(f) \colon \m{BSC}_k(X) \to k$ of $f$. Moreover, taking a representative $x_n \in \Psi(n)$ for each $n \in \N$, one obtains a continuous map $x \colon \N \to X \hookrightarrow \m{BSC}_k(X)$. Since $\m{BSC}_k(X)$ is an object of $X/\m{TDCHTop}$, a unique continuous extension $\m{BSC}_k(x) \colon \m{UF}(\N) \cong \m{BSC}_k(\N) \to \m{BSC}_k(X)$ of $x$ exists. The composite $\m{BSC}_k(f) \circ \m{BSC}_k(x) \colon \m{UF}(\N) \to \m{BSC}_k(X) \to k$ is a continuous extension of the composite $f \circ x = \varphi \circ \psi \colon \N \to k$. In particular, $\m{BSC}_k(f) \circ \m{BSC}_k(x)$ is continuous at $\mathscr{F} \in \m{UF}(\N)$, but it contradicts the fact that $\varphi \circ \psi$ is an injective map whose image is a closed discrete subspace. An injective net whose image is discrete and closed never has a limit. It is a contradiction. Therefore one obtains $\F = \Q$, and (ii) holds by the assumption. Take a $k$-linear combination $g = \sum_{i = 1}^{n} a_i1_{U_i} \in \m{C}_{\m{bd}}(X,k)$ of idempotents with $\| f - g \| < 1$. Now the image of $g$ contains at most $n$ points, and hence there is an integer $m \in \N$ such that $g(x) \notin \psi(m)$ for any $x \in X$ identifying the cosets $\tilde{k}$ as a family of disjoint clopen subsets of $k^{\circ}$ in the tautological sense. Then one has $|f(x_m) - g(x_m)| = 1$, and it contradicts the condition $\| f - g \| < 1$. Thus $X$ is compact, or $k$ is a local field or a finite field.
\end{proof}

Since we did not use the assumption that $\F \neq \Q$ in the proof that (ii) with $\F \neq \Q$ implies (iii), the condition that $\iota_{k/\F}$ is an isometric isomorphism is weaker than (ii).

\subsection{Non-Archimedean Gel'fand Theory}
\label{Non-Archimedean Gel'fand Theory}

We establish non-Archimedean Gel'fand theory for a zero-dimensional Hausdorff space. We recall that a completely regular Hausdorff space is a topological space which can be embedded in a direct product of copies of the closed unit disc $\C^{\circ} \subset \C$ as a subspace. On the other hand, a non-Archimedean counterpart of a completely regular Hausdorff space over $k$ is a topological space which can be embedded in a direct product of copies of the closed unit disc $k^{\circ} \subset k$ as a subspace. We call such a topological space a {\it non-Archimedean completely regular Hausdorff space} over $k$. A direct product of copies of $k^{\circ}$ is a zero-dimensional Hausdorff space. A subset of a zero-dimensional Hausdorff space is again a zero-dimensional Hausdorff space, and so is a non-Archimedean completely regular Hausdorff space over $k$. Now we verify that the converse also holds in the case where $k$ is a local field or a finite field.

\begin{lmm}
\label{completely regular Hausdorff}
Suppose that $k$ is a local field or a finite field endowed with the trivial norm. The following are equivalent:
\begin{itemize}
\item[(i)] The space $X$ is zero-dimensional and Hausdorff.
\item[(ii)] The space $X$ is Hausdorff, and bounded continuous $k$-valued functions separates a point and a disjoint closed subset of $X$, i.e.\ for any $x \in X$ and any closed subset $F \subset X$ with $x \notin F$, there is an $f \in \m{C}_{\m{bd}}(X,k)$ such that $f(x) = 0$ and $f(y) = 1$ for any $y \in F$.
\item[(iii)] The continuous map $\iota'_k \colon X \to \m{SC}_k(X)$ is a homeomorphism onto the image.
\item[(iv)] The space $X$ is embedded in a direct product of copies of $k^{\circ}$ as a subspace.
\end{itemize}
\end{lmm}

We note that the description of (ii)-(iv) seem to depend on the base field $k$ while (i) does not. Therefore the notion of ``a non-Archimedean completely regular Hausdorff space'' is independent of the base field.

\begin{proof}
Recall that $\m{SC}_k(X)$ is a closed subspace of a direct product of copies of $k^{\circ}$, and hence (iii) implies (iv). Moreover, (iii) implies (i) as we mentioned in the beginning of this section.

Assume (i). We verify (ii). Let $x \in X$ and $F \subset X$ be a closed subset with $x \notin F$. Since $X$ is zero-dimensional, there is a clopen neighbourhood $U \subset X$ of $x$ contained in the open subset $X \backslash F \subset X$, and the characteristic function $1_{X \backslash U}$ separates $x$ and $F$.

Assume (ii). We verify (iii). Since $X$ is Hausdorff, a point of $X$ is closed. For $x,y \in X$ with $x \neq y$, take an $f \in \m{C}_{\m{bd}}(X,k)$ which separates $x$ and $y$. Then $f \neq 0$. Since the valuation of $k$ is discrete or trivial, the image $|k| \subset \lbrack 0,\infty )$ is closed. By the definition of the supremum norm, $\| \m{C}_{\m{bd}}(X,k) \| \subset [0,\infty)$ is contained in the closure of $|k| \subset [0,\infty)$, and hence $\| \m{C}_{\m{bd}}(X,k) \| \subset |k|$. Therefore there is an $a \in k^{\times}$ such that $0 \neq \| f \| = |a|$. Then one has $\| a^{-1}f \| = 1$ and $a^{-1}f \in \m{C}_{\m{bd}}(X,k)(1)$ separates $x$ and $y$. Therefore one has $\iota'_k(x) \neq \iota'_k(y)$ comparing their $(a^{-1}f)$-th entry, and $\iota'_k$ is injective. In order to prove that $\iota'_k$ is an open map onto the image, take an open subset $U \subset X$. For an $x \in U$, take an $f \in \m{C}_{\m{bd}}(X,k)$ such that $f(x) = 0$ and $f(y) = 1$ for any $y \in X \backslash U$. By the same argument as above, there is an $a \in k^{\times}$ such that $\| f \| = |a|$. Then the pre-image by $\iota'_k$ of the open subset $V \subset \m{SC}_k(X)$ given by the condition that the $(a^{-1}f)$-th entry is contained in the open neighbourhood $k \backslash \{ a^{-1} \} \subset k$ of $0 \in k$ is an open neighbourhood of $x$ contained in $U$. Therefore the image $\iota'_k(U)$ contains the open neighbourhood $V \cap \iota'_k(X)$ of $\iota'_k(x)$ in $\iota'_k(X)$, and thus $\iota'_k(U) \subset \iota'_k(X)$ is open. We conclude that $\iota'_k$ is a homeomorphism onto the image.
\end{proof}

\begin{prp}
\label{evaluation map}
The map $\iota_k \colon X \to \m{BSC}_k(X)$ is a homeomorphism onto the image if and only if $X$ is zero-dimensional and Hausdorff.
\end{prp}

\begin{proof}
By Proposition \ref{independent of k}, it is reduced to the case where $k = \Q_p$ for a prime number $p \in \N$. The assertion immediately follows from Corollary \ref{SC and BSC} (i) and Lemma \ref{completely regular Hausdorff}.
\end{proof}

\begin{dfn}
Let $\A \subset \m{C}_{\m{bd}}(X,k)$ be a closed $k$-subalgebra. For $x,x' \in X$, we write $x \sim_{\A} x'$ if $f(x) = f(x')$ for any $f \in \A$. The binary relation $\sim_{\A}$ is an equivalence relation, and we denote by $X/\A$ the quotient space $X/\sim_{\A}$. We say $\A$ separates points of $X$ if the condition $x \sim_{\A} x'$ implies $x= x'$ for any $x,x' \in X$.
\end{dfn}

\begin{lmm}
\label{maximal separability}
The map $\iota'_k \colon X \to \m{SC}_k(X)$ uniquely factors through the canonical projection $X \twoheadrightarrow X/\m{C}_{\m{bd}}(X,k)$, and the induced map $X/\m{C}_{\m{bd}}(X,k) \to \m{SC}_k(X)$ is an injective continuous map.
\end{lmm}

\begin{proof}
It immediately follows from the definitions of $\sim_{\m{C}_{\m{bd}}(X,k)}$ and $\m{SC}_k(X)$.
\end{proof}

\begin{lmm}
Suppose that $X$ is zero-dimensional and Hausdorff. Then $\m{C}_{\m{bd}}(X,k)$ separates points of $X$.
\end{lmm}

\begin{proof}
By Lemma \ref{completely regular Hausdorff} and Lemma \ref{maximal separability}, the projection $X \twoheadrightarrow X/\m{C}_{\m{bd}}(X,k)$ is injective.
\end{proof}

\begin{dfn}
\label{faithful}
A topological space $(Y,f)$ under $X$ is said to be {\it faithful} if $f \colon X \to Y$ is a homeomorphism onto the image, to be {\it full} if $f(X)$ is dense in $Y$, and to be {\it fully faithful} if it is full and faithful. We denote by $X/\m{TDCHTop}_{\m{ff}} \subset X/\m{TDCHTop}$ the full subcategory of fully faithful totally disconnected compact Hausdorff spaces under $X$.
\end{dfn}

By Lemma \ref{completely regular Hausdorff}, $X/\m{TDCHTop}_{\m{ff}}$ is a non-empty category if and only if $X$ is zero-dimensional and Hausdorff. Suppose that $X$ is zero-dimensional and Hausdorff. The isomorphism relation in $X/\m{TDCHTop}_{\m{ff}}$ is an equivalence relation in a class. We denote by $\mathscr{C}(X)$ the class $(X/\m{TDCHTop}_{\m{ff}})/\cong$ of equivalence classes. The class $\mathscr{C}(X)$ is not a proper class. Indeed, for any $(Y,f) \in \m{ob}(X/\m{TDCHTop}_{\m{ff}})$, $f$ extends to $\tilde{f} \colon \m{UF}(X) \to Y$ by Theorem \ref{model}. Since $\m{UF}(X)$ is compact and $Y$ is Hausdorff, $\tilde{f}(\m{UF}(X)) \subset Y$ is a closed subspace containing the dense subspace $f(X) \subset Y$, and hence $\tilde{f}$ is a surjective closed map. Therefore $(Y,f)$ is obtained as a quotient of $\m{UF}(X)$, and $\mathscr{C}(X)$ admits a set-theoretical representative.

\begin{thm}
\label{gelfand}
Suppose that $k$ is a local field or a finite field endowed with the trivial norm, and that $X$ is zero-dimensional and Hausdorff. Then there is a contravariant-functorial one-to-one correspondence between $\mathscr{C}(X)$ and the set of closed $k$-subalgebras of $\m{C}_{\m{bd}}(X,k)$ separating points of $X$.
\end{thm}

\begin{proof}
Denote by $\mathscr{C}'(X)$ the set of closed $k$-subalgebras of $\m{C}_{\m{bd}}(X,k)$ separating points of $X$. The correspondences are given in the following way:
\begin{eqnarray*}
\begin{array}{ccc}
  \mathscr{C}(X) & \longleftrightarrow & \mathscr{C}'(X)\\
  \left[ f \colon X \hookrightarrow Y \right] & \rightsquigarrow & \m{Im}(f^{\m{a}} \colon \m{C}(Y,k) \hookrightarrow \m{C}_{\m{bd}}(X,k))\\
  \left[ X \hookrightarrow \m{BSC}_k(X) \twoheadrightarrow \M_k(\A) \right] & \leftsquigarrow & (\A \subset \m{C}_{\m{bd}}(X,k)).
\end{array}
\end{eqnarray*}
They are the inverses of each other by the generalised Stone--Weierstrass theorem (\cite{Ber1} 9.2.5.\ Theorem). We remark that for any fully faithful totdally disconnected compact Hausdorff space $(Y,f)$ under $X$, the associated bounded homomorphism $f^{\m{a}} \colon \m{C}(Y,k) \hookrightarrow \m{C}_{\m{bd}}(X,k)$ is an isometry because $X$ is a dense subspace of $Y$, and hence its image is closed. On the other hand for any closed $k$-subalgebra $\A$ of $\m{C}_{\m{bd}}(X,k)$ separating points of $X$, the associated continuous map $X \hookrightarrow \m{BSC}_k(X) \twoheadrightarrow \M_k(\A)$ is a homeomorphism onto the image which is a dense subspace, because $X$ is a dense subspace of $\m{BSC}_k(X)$, the condition that $\A$ separates points of$X$ guarantees the injectivity, and every continuous map between compact Hausdorff spaces $\m{BSC}_k(X)$ and $\M_k(\A)$ is a closed map.
\end{proof}

\vspace{0.4in}
\addcontentsline{toc}{section}{Acknowledgements}
\noindent {\Large \bf Acknowledgements}
\vspace{0.1in}

I would like to thank Takeshi Tsuji. I appreciate having gracious discussions in seminars. I am profoundly grateful to Atsushi Matsuo for their great helps and instructive advices on writing this paper. I express my gratitude to Atsushi Yamashita for joining in seminars on this paper. I am thankful to my friends for devoting time to frequent discourse. I acknowledge my family's deep affection.

I am a research fellow of Japan Society for the Promotion of Science. This work was supported by the Program for Leading Graduate 
Schools, MEXT, Japan.


\addcontentsline{toc}{section}{References}
\begin{thebibliography}{99}

\bibitem[Ban]{Ban} Bernhard Banaschewski, {\it \"Uber Nulldimensionale R\"aume}, Mathematische Nachrichten, Volume 13, p.\ 129-140, 1955.

\bibitem[Ber1]{Ber1} Vladimir G.\ Berkovich, {\it Spectral Theory and Analytic Geometry over non-Archimedean Fields}, Mathematical Surveys and Monographs, Number 33, the American Mathematical Society, 1990.

\bibitem[Ber2]{Ber2} Vladimir G.\ Berkovich, {\it \'Etale Cohomology for non-Archimedean Analytic Spaces}, Publications Math\'ematiques de l'Institut des Hautes \'Etudes Scientifiques, Number 1, p.\ 5-161, Springer, 1993.

\bibitem[BGR]{BGR} S.\ Bosch, U.\ G\"untzer, and R.\ Remmert, {\it Non-Archimedean Analysis A Systematic Approach to Rigid Analytic Geometry}, Grundlehren der mathematischen Wissenschaften 261, A Series of Comprehensive Studies in Mathematics, Springer, 1984.

\bibitem[BJ]{BJ} Francis Borceux, George Janelidze, {\it Galois Theories}, Cambridge Studies in Advanced Mathematics 72, Cambridge University Press, 2001.

\bibitem[Bou]{Bou} Nicolas Bourbaki, {\it Espaces Vectoriels Topologiques}, El\'ements de math\'ematique, Volume V, Hermann, 1953.

\bibitem[Dal]{Dal} H.\ Garth Dales {\it Banach Algebra and Automatic Continuity}, London Mathematical Society Monographs New Series, Volume 24, Oxford Science Publications, 2000.

\bibitem[Die]{Die} Jean Dieudonn\'e, {\it Sur les Fonctions Continues $p$-adiques}, Bulletin des Sciences Mathematiques, Volume 68, p.\ 79-95, 1944.

\bibitem[Dou]{Dou} Ronald G.\ Douglas, {\it Banach Algebra Techniques in Operator Theory}, Pure and Applied Mathematics, Volume 49, Academic Press, 1972.

\bibitem[Eng]{Eng} Ryszard Engelking, {\it General Topology}, Pa\'nstwowe Wydawnictwo Naukowe, 1977.

\bibitem[EM1]{EM1} Alain Escassut, Nicolas Mainetti, {\it Multiplicative Spectrum of Ultrametric Banach Algebras of Continuous Functions}, Topology and Its Applications, Volume 157, p.\ 2505-2515, 2010.

\bibitem[EM2]{EM2} Alain Escassut, Nicolas Mainetti, {\it Morphisms between ultrametric Banach algebras and maximal ideals of finite codimension}, Comtemporary Mathematics, Advances in Ultrametric Analysis, Volume 596, p.\ 63-71, 2013.

\bibitem[Gue]{Gue} Bernard Guennebaud, {\it th\`ese Universit\'e de Poitiers}, Sur une Notion de Spectre pour les Alg\`ebres Norm\'ees Ultram\'etriques, 1973.

\bibitem[Joh]{Joh} Peter T.\ Johnstone, {\it Stone Spaces}, Cambridge Studies in Advanced Mathematics 3, Cambridge University Press, 1982.

\bibitem[Kap]{Kap} Irving Kaplansky, {\it The Weierstrass Theorem in Fields with Valuations}, Proceedings of the American Mathematical Society, Volume 1, p.\ 356-357, 1950.

\bibitem[Kur]{Kur} Kazimier Kuratowski, {\it Topologie Volume I}, Pa\'nstwowe Wydawnictwo Naukowe, 1958.

\bibitem[Par]{Par} Parfeny P.\ Saworotnow, {\it Totally Disconnected Compactification}, International Journal of Mathematics and Mathematical Sciences, Volume 16, Number 4, p.\ 653-656, 1993.

\bibitem[Sol]{Sol} R.\ M.\ Solovay, {\it Discontinuous Homomorphisms of Banach Algebras}, preprint, 1976.

\bibitem[Sto1]{Sto1} Arthur H.\ Stone, {\it Paracompactness and Product Spaces}, Bulletin of the American Mathematical Society, Volume 54, p.\ 977-982, 1948.

\bibitem[Sto2]{Sto2} Marshall H.\ Stone, {\it Boolean Algebras and Their Application to Topology}, Proceedings of National Academy of Sciences of the United State of America, Volume 20, Number 3, p.\ 197-202, 1934.

\bibitem[Sto3]{Sto3} Marshall H.\ Stone, {\it The Theory of Representations of Boolean Algebras}, Transactions of the American Mathematical Society, Volume 40, Number 1, p.\ 37-111, 1936.

\bibitem[Tar]{Tar} Abolfazl Tarizadeh, {\it On the Category of Profinite Spaces as a Reflective Subcategory}, preprint.

\bibitem[Woo]{Woo} W.\ H.\ Woodin, {\it A Discontinuous Homomorphism from $C(X)$ without CH}, Journal of the London Mathematical Society, Volume s2-48 (2), p.\ 299-315, 1993.

\end{thebibliography}
\end{document}